\documentclass[11pt]{article}
\usepackage[margin=1in]{geometry}
\usepackage{algorithm}
\usepackage{algpseudocode}
\usepackage{url} 

\providecommand{\headers}[2]{} 

\let\oldtitle\title
\renewcommand{\title}[2][]{\oldtitle{#2}}

\makeatletter
\newcommand{\LoadPackageUnlessAcm}[2][]{%
  \@ifclassloaded{acmart}{%
    \PackageInfo{preamble}{Skipping package `#2' because acmart class is used}%
  }{%
    \if\relax\detokenize{#1}\relax
      \usepackage{#2}%
    \else
      \usepackage[#1]{#2}%
    \fi
  }%
}
\makeatother

\LoadPackageUnlessAcm{caption}
\LoadPackageUnlessAcm{subcaption}
\LoadPackageUnlessAcm{authblk}
\LoadPackageUnlessAcm{hyperref} 

\usepackage{amstext}
\usepackage{stmaryrd}
\usepackage{algpseudocode}

\usepackage{units}
\usepackage{cancel}
\usepackage{graphicx}

\usepackage{xcolor}
\usepackage{array}
\usepackage{verbatim}
\usepackage{booktabs}
\usepackage{calc}
\usepackage{textcomp}
\usepackage{multirow}
\usepackage{tikz,tikzscale}

\usepackage{pgfplots}
\usepackage{pgfplotstable}
\usepackage{pgfkeys}
\pgfplotsset{compat=1.18}

\definecolor{deepgreen}{RGB}{0,100,0}

\newif\ifrunscripts
\runscriptstrue
\newcommand{\executeextract}[3]{%
    \ifrunscripts
        \immediate\write18{python3 scripts/extract_data.py --input "#1" --query "#2" --output "#3"}%
    \fi
}

\executeextract{results/struct_2d_SRichardson/summary_braess.csv}{SELECT distortion, avg FROM data WHERE degree=2 AND mu=1.0 AND schur=1 AND smooth=1}{results/tmp/struct_2d_braess_p2.csv}
\executeextract{results/struct_2d_SRichardson/summary_braess.csv}{SELECT distortion, avg FROM data WHERE degree=3 AND mu=1.0 AND schur=1 AND smooth=1}{results/tmp/struct_2d_braess_p3.csv}
\executeextract{results/struct_2d_SRichardson/summary_braess.csv}{SELECT distortion, avg FROM data WHERE degree=7 AND mu=1.0 AND schur=1 AND smooth=1}{results/tmp/struct_2d_braess_p7.csv}

\executeextract{results/struct_2d_SRichardson/summary_block.csv}{SELECT distortion, avg FROM data WHERE degree=2 AND mu=1.0 AND schur=1 AND smooth=1}{results/tmp/struct_2d_block_p2.csv}
\executeextract{results/struct_2d_SRichardson/summary_block.csv}{SELECT distortion, avg FROM data WHERE degree=3 AND mu=1.0 AND schur=1 AND smooth=1}{results/tmp/struct_2d_block_p3.csv}
\executeextract{results/struct_2d_SRichardson/summary_block.csv}{SELECT distortion, avg FROM data WHERE degree=7 AND mu=1.0 AND schur=1 AND smooth=1}{results/tmp/struct_2d_block_p7.csv}

\executeextract{results/struct_3d_SRichardson/summary_block.csv}{SELECT distortion, avg FROM data WHERE degree=2 AND mu=1.0 AND schur=1 AND smooth=1}{results/tmp/struct_3d_block_p2.csv}
\executeextract{results/struct_3d_SRichardson/summary_block.csv}{SELECT distortion, avg FROM data WHERE degree=3 AND mu=1.0 AND schur=1 AND smooth=1}{results/tmp/struct_3d_block_p3.csv}
\executeextract{results/struct_3d_SRichardson/summary_block.csv}{SELECT distortion, avg FROM data WHERE degree=7 AND mu=1.0 AND schur=1 AND smooth=1}{results/tmp/struct_3d_block_p7.csv}

\executeextract{results/struct_2d_exact/summary_block.csv}{SELECT distortion, avg FROM data WHERE degree=2 AND mu=1.0 AND schur=1 AND smooth=1}{results/tmp/struct_2d_exact_p2.csv}
\executeextract{results/struct_2d_exact/summary_block.csv}{SELECT distortion, avg FROM data WHERE degree=3 AND mu=1.0 AND schur=1 AND smooth=1}{results/tmp/struct_2d_exact_p3.csv}
\executeextract{results/struct_2d_exact/summary_block.csv}{SELECT distortion, avg FROM data WHERE degree=7 AND mu=1.0 AND schur=1 AND smooth=1}{results/tmp/struct_2d_exact_p7.csv}

\executeextract{results/struct_3d_exact/summary_block.csv}{SELECT distortion, avg FROM data WHERE degree=2 AND mu=1.0 AND schur=1 AND smooth=1}{results/tmp/struct_3d_exact_p2.csv}
\executeextract{results/struct_3d_exact/summary_block.csv}{SELECT distortion, avg FROM data WHERE degree=3 AND mu=1.0 AND schur=1 AND smooth=1}{results/tmp/struct_3d_exact_p3.csv}
\executeextract{results/struct_3d_exact/summary_block.csv}{SELECT distortion, avg FROM data WHERE degree=7 AND mu=1.0 AND schur=1 AND smooth=1}{results/tmp/struct_3d_exact_p7.csv}

\executeextract{results/simplex_2d_SRichardson/summary_braess.csv}{SELECT distortion, avg FROM data WHERE degree=2 AND mu=1.0 AND schur=1 AND smooth=1}{results/tmp/simplex_2d_braess_p2.csv}
\executeextract{results/simplex_2d_SRichardson/summary_braess.csv}{SELECT distortion, avg FROM data WHERE degree=3 AND mu=1.0 AND schur=1 AND smooth=1}{results/tmp/simplex_2d_braess_p3.csv}
\executeextract{results/simplex_2d_SRichardson/summary_braess.csv}{SELECT distortion, avg FROM data WHERE degree=7 AND mu=1.0 AND schur=1 AND smooth=1}{results/tmp/simplex_2d_braess_p7.csv}

\executeextract{results/simplex_2d_SRichardson/summary_block.csv}{SELECT distortion, avg FROM data WHERE degree=2 AND mu=1.0 AND schur=1 AND smooth=1}{results/tmp/simplex_2d_block_p2.csv}
\executeextract{results/simplex_2d_SRichardson/summary_block.csv}{SELECT distortion, avg FROM data WHERE degree=3 AND mu=1.0 AND schur=1 AND smooth=1}{results/tmp/simplex_2d_block_p3.csv}
\executeextract{results/simplex_2d_SRichardson/summary_block.csv}{SELECT distortion, avg FROM data WHERE degree=7 AND mu=1.0 AND schur=1 AND smooth=1}{results/tmp/simplex_2d_block_p7.csv}

\executeextract{results/simplex_3d_SRichardson/summary_block.csv}{SELECT distortion, avg FROM data WHERE degree=2 AND mu=1.0 AND schur=1 AND smooth=1}{results/tmp/simplex_3d_block_p2.csv}
\executeextract{results/simplex_3d_SRichardson/summary_block.csv}{SELECT distortion, avg FROM data WHERE degree=3 AND mu=1.0 AND schur=1 AND smooth=1}{results/tmp/simplex_3d_block_p3.csv}
\executeextract{results/simplex_3d_SRichardson/summary_block.csv}{SELECT distortion, avg FROM data WHERE degree=7 AND mu=1.0 AND schur=1 AND smooth=1}{results/tmp/simplex_3d_block_p7.csv}

\executeextract{results/simplex_3d_SRichardson/summary_braess.csv}{SELECT distortion, avg FROM data WHERE degree=2 AND mu=1.0 AND schur=1 AND smooth=1}{results/tmp/simplex_3d_braess_p2.csv}
\executeextract{results/simplex_3d_SRichardson/summary_braess.csv}{SELECT distortion, avg FROM data WHERE degree=3 AND mu=1.0 AND schur=1 AND smooth=1}{results/tmp/simplex_3d_braess_p3.csv}
\executeextract{results/simplex_3d_SRichardson/summary_braess.csv}{SELECT distortion, avg FROM data WHERE degree=7 AND mu=1.0 AND schur=1 AND smooth=1}{results/tmp/simplex_3d_braess_p7.csv}

\executeextract{results/simplex_3d_SRichardson/summary_block.csv}{SELECT distortion, avg FROM data WHERE degree=2 AND mu=1.0 AND schur=1 AND smooth=1}{results/tmp/simplex_3d_block_p2.csv}
\executeextract{results/simplex_3d_SRichardson/summary_block.csv}{SELECT distortion, avg FROM data WHERE degree=3 AND mu=1.0 AND schur=1 AND smooth=1}{results/tmp/simplex_3d_block_p3.csv}
\executeextract{results/simplex_3d_SRichardson/summary_block.csv}{SELECT distortion, avg FROM data WHERE degree=7 AND mu=1.0 AND schur=1 AND smooth=1}{results/tmp/simplex_3d_block_p7.csv}

\executeextract{results/struct_2d_SRichardson/summary_braess.csv}{SELECT mu, avg FROM data WHERE degree=2 AND distortion=0.0 AND schur=1 AND smooth=1}{results/tmp/struct_2d_braess_jump_p2.csv}
\executeextract{results/struct_2d_SRichardson/summary_braess.csv}{SELECT mu, avg FROM data WHERE degree=3 AND distortion=0.0 AND schur=1 AND smooth=1}{results/tmp/struct_2d_braess_jump_p3.csv}
\executeextract{results/struct_2d_SRichardson/summary_braess.csv}{SELECT mu, avg FROM data WHERE degree=7 AND distortion=0.0 AND schur=1 AND smooth=1}{results/tmp/struct_2d_braess_jump_p7.csv}

\executeextract{results/simplex_2d_SRichardson/summary_braess.csv}{SELECT mu, avg FROM data WHERE degree=2 AND distortion=0.0 AND schur=1 AND smooth=1}{results/tmp/simplex_2d_braess_jump_p2.csv}
\executeextract{results/simplex_2d_SRichardson/summary_braess.csv}{SELECT mu, avg FROM data WHERE degree=3 AND distortion=0.0 AND schur=1 AND smooth=1}{results/tmp/simplex_2d_braess_jump_p3.csv}
\executeextract{results/simplex_2d_SRichardson/summary_braess.csv}{SELECT mu, avg FROM data WHERE degree=7 AND distortion=0.0 AND schur=1 AND smooth=1}{results/tmp/simplex_2d_braess_jump_p7.csv}

\executeextract{results/struct_2d_SRichardson/summary_block.csv}{SELECT mu, avg FROM data WHERE degree=2 AND distortion=0.0 AND schur=1 AND smooth=1}{results/tmp/struct_2d_block_jump_p2.csv}
\executeextract{results/struct_2d_SRichardson/summary_block.csv}{SELECT mu, avg FROM data WHERE degree=3 AND distortion=0.0 AND schur=1 AND smooth=1}{results/tmp/struct_2d_block_jump_p3.csv}
\executeextract{results/struct_2d_SRichardson/summary_block.csv}{SELECT mu, avg FROM data WHERE degree=7 AND distortion=0.0 AND schur=1 AND smooth=1}{results/tmp/struct_2d_block_jump_p7.csv}

\executeextract{results/struct_2d_SRichardson/summary_braess.csv}{SELECT distortion, avg FROM data WHERE degree=2 AND mu=10000.0 AND schur=1 AND smooth=1}{results/tmp/struct_2d_braess_dist_mu4_p2.csv}
\executeextract{results/struct_2d_SRichardson/summary_braess.csv}{SELECT distortion, avg FROM data WHERE degree=3 AND mu=10000.0 AND schur=1 AND smooth=1}{results/tmp/struct_2d_braess_dist_mu4_p3.csv}
\executeextract{results/struct_2d_SRichardson/summary_braess.csv}{SELECT distortion, avg FROM data WHERE degree=7 AND mu=10000.0 AND schur=1 AND smooth=1}{results/tmp/struct_2d_braess_dist_mu4_p7.csv}

\executeextract{results/simplex_2d_SRichardson/summary_braess.csv}{SELECT distortion, avg FROM data WHERE degree=2 AND mu=10000.0 AND schur=1 AND smooth=1}{results/tmp/simplex_2d_braess_dist_mu4_p2.csv}
\executeextract{results/simplex_2d_SRichardson/summary_braess.csv}{SELECT distortion, avg FROM data WHERE degree=3 AND mu=10000.0 AND schur=1 AND smooth=1}{results/tmp/simplex_2d_braess_dist_mu4_p3.csv}
\executeextract{results/simplex_2d_SRichardson/summary_braess.csv}{SELECT distortion, avg FROM data WHERE degree=7 AND mu=10000.0 AND schur=1 AND smooth=1}{results/tmp/simplex_2d_braess_dist_mu4_p7.csv}

\executeextract{results/struct_2d_SRichardson/summary_block.csv}{SELECT distortion, avg FROM data WHERE degree=2 AND mu=10000.0 AND schur=1 AND smooth=1}{results/tmp/struct_2d_block_dist_mu4_p2.csv}
\executeextract{results/struct_2d_SRichardson/summary_block.csv}{SELECT distortion, avg FROM data WHERE degree=3 AND mu=10000.0 AND schur=1 AND smooth=1}{results/tmp/struct_2d_block_dist_mu4_p3.csv}
\executeextract{results/struct_2d_SRichardson/summary_block.csv}{SELECT distortion, avg FROM data WHERE degree=7 AND mu=10000.0 AND schur=1 AND smooth=1}{results/tmp/struct_2d_block_dist_mu4_p7.csv}

\executeextract{results/struct_2d_SRich_nonstat/summary_braess.csv}{SELECT mu, avg FROM data WHERE degree=2 AND distortion=0.0 AND schur=1 AND smooth=1}{results/tmp/struct_2d_braess_jump_nonstat_p2.csv}
\executeextract{results/struct_2d_SRich_nonstat/summary_braess.csv}{SELECT mu, avg FROM data WHERE degree=3 AND distortion=0.0 AND schur=1 AND smooth=1}{results/tmp/struct_2d_braess_jump_nonstat_p3.csv}
\executeextract{results/struct_2d_SRich_nonstat/summary_braess.csv}{SELECT mu, avg FROM data WHERE degree=7 AND distortion=0.0 AND schur=1 AND smooth=1}{results/tmp/struct_2d_braess_jump_nonstat_p7.csv}

\executeextract{results/struct_2d_SRich_nonstat/summary_block.csv}{SELECT mu, avg FROM data WHERE degree=2 AND distortion=0.0 AND schur=1 AND smooth=1}{results/tmp/struct_2d_block_jump_nonstat_p2.csv}
\executeextract{results/struct_2d_SRich_nonstat/summary_block.csv}{SELECT mu, avg FROM data WHERE degree=3 AND distortion=0.0 AND schur=1 AND smooth=1}{results/tmp/struct_2d_block_jump_nonstat_p3.csv}
\executeextract{results/struct_2d_SRich_nonstat/summary_block.csv}{SELECT mu, avg FROM data WHERE degree=7 AND distortion=0.0 AND schur=1 AND smooth=1}{results/tmp/struct_2d_block_jump_nonstat_p7.csv}

\executeextract{results/struct_2d_SRich_nonstat/summary_braess.csv}{SELECT distortion, avg FROM data WHERE degree=2 AND mu=10000.0 AND schur=1 AND smooth=1}{results/tmp/struct_2d_braess_dist_nonstat_mu4_p2.csv}
\executeextract{results/struct_2d_SRich_nonstat/summary_braess.csv}{SELECT distortion, avg FROM data WHERE degree=3 AND mu=10000.0 AND schur=1 AND smooth=1}{results/tmp/struct_2d_braess_dist_nonstat_mu4_p3.csv}
\executeextract{results/struct_2d_SRich_nonstat/summary_braess.csv}{SELECT distortion, avg FROM data WHERE degree=7 AND mu=10000.0 AND schur=1 AND smooth=1}{results/tmp/struct_2d_braess_dist_nonstat_mu4_p7.csv}

\executeextract{results/struct_2d_SRich_nonstat/summary_block.csv}{SELECT distortion, avg FROM data WHERE degree=2 AND mu=10000.0 AND schur=1 AND smooth=1}{results/tmp/struct_2d_block_dist_nonstat_mu4_p2.csv}
\executeextract{results/struct_2d_SRich_nonstat/summary_block.csv}{SELECT distortion, avg FROM data WHERE degree=3 AND mu=10000.0 AND schur=1 AND smooth=1}{results/tmp/struct_2d_block_dist_nonstat_mu4_p3.csv}
\executeextract{results/struct_2d_SRich_nonstat/summary_block.csv}{SELECT distortion, avg FROM data WHERE degree=7 AND mu=10000.0 AND schur=1 AND smooth=1}{results/tmp/struct_2d_block_dist_nonstat_mu4_p7.csv}

\executeextract{results/struct_3d_SRichardson/summary_braess.csv}{SELECT distortion, avg FROM data WHERE degree=2 AND mu=1.0 AND schur=1 AND smooth=1}{results/tmp/struct_3d_braess_p2.csv}
\executeextract{results/struct_3d_SRichardson/summary_braess.csv}{SELECT distortion, avg FROM data WHERE degree=3 AND mu=1.0 AND schur=1 AND smooth=1}{results/tmp/struct_3d_braess_p3.csv}
\executeextract{results/struct_3d_SRichardson/summary_braess.csv}{SELECT distortion, avg FROM data WHERE degree=7 AND mu=1.0 AND schur=1 AND smooth=1}{results/tmp/struct_3d_braess_p7.csv}

\executeextract{results/poisson_Ref/Jac_cart_2D.csv}{SELECT distortion, avg FROM data WHERE degree=3}{results/tmp/poisson_2d_p3.csv}
\executeextract{results/poisson_Ref/Jac_cart_2D.csv}{SELECT distortion, avg FROM data WHERE degree=7}{results/tmp/poisson_2d_p7.csv}

\pgfplotsset{
  p2/.style={teal, mark=*},
  p3/.style={orange, mark=triangle*},
  p7/.style={blue, mark=square*},
}

\newcommand{\plotwidth}{0.49\columnwidth}
\newcommand{\plotheight}{0.4\columnwidth}

\pgfplotscreateplotcyclelist{paper markers}{
  {teal,    solid, mark=*},
  {orange,   solid, mark=triangle*},
  {blue,     solid, mark=square*},
  {red,      solid, mark=pentagon*},
  {violet,     solid, mark=otimes*},
  {brown,    solid, mark=halfcircle*},
  {black,   solid, mark=diamond*},
}

\colorlet{BarColorOne}{orange!60!white}
\colorlet{BarColorTwo}{teal!80!white}
\colorlet{BarColorThree}{red!80!black}
\colorlet{BarColorFour}{blue!60!white}
\colorlet{BarColorFive}{brown!70!black}

\pgfplotsset{legend swatch/.style={area legend, draw=none}}

\pgfplotsset{
  paperplot/.style={
      width=\plotwidth,
      height=\plotheight,
      cycle list name=paper markers,
      every axis/.append style={font=\small},
      grid=major,
      thick
    }
}

\pgfplotsset{
  longplot/.style={
      width=0.8\columnwidth,
      height=\plotheight,
      cycle list name=paper markers,
      every axis/.append style={font=\small},
      grid=major,
      thick
    }
}

\usepackage[all]{xy}

\usepackage{mathtools}
\usepackage{placeins}

\usetikzlibrary{calc}
\usetikzlibrary{shadings}

\usepgfplotslibrary{fillbetween}

\usepackage[mode=multiuser,status=draft]{fixme} 
\fxsetup{innerlayout=inline}
\fxsetup{targetlayout=colorcb}

\fxusetheme{color}
\FXRegisterAuthor{mw}{envmw}{MW}

\definecolor{apply}{rgb}{0.3,.7,0.}
\definecolor{applyface}{rgb}{.9,.95,0.}
\definecolor{invert}{rgb} {0.5,0.,1.}
\colorlet{invertface}{invert!60!red}

\usepackage[colorinlistoftodos,prependcaption,textsize=small]{todonotes}

\def\mesh{\mathbb M}

\def\R{\mathbb R}

\def\resth#1{I_{#1}^{\downarrow_{h}}}
\def\prolh#1{I_{#1}^{\uparrow_{h}}}

\def\prolp#1{I_{#1}^{\Uparrow_{p}}}
\def\restp#1{I_{#1}^{\Downarrow_{p}}}

\usepackage{colortbl} 
\newcommand{\myrowcolor}{\rowcolor[gray]{0.925}}

\usepackage{booktabs} 
\usepackage[font=footnotesize,position=b]{subcaption}
\usepackage{adjustbox}

\newcommand{\pluseq}{\stackrel{+}{\gets}}

\newcommand{\LUpNew}{local\_update}

\newcommand{\residual}{r}
\newcommand{\res}{\text{res}}

\makeatletter
\def\pgfplots@getautoplotspec into#1{%
    \begingroup
    \let#1=\pgfutil@empty
    \pgfkeysgetvalue{/pgfplots/cycle multi list/@dim}\pgfplots@cycle@dim
    \let\pgfplots@listindex=\pgfplots@numplots
    \pgfkeysgetvalue{/pgfplots/cycle list set}\pgfplots@listindex@set
    \ifx\pgfplots@listindex@set\pgfutil@empty
    \else
      \c@pgf@counta=\pgfplots@listindex
      \c@pgf@countb=\pgfplots@listindex@set
      \advance\c@pgf@countb by -\c@pgf@counta
      \globaldefs=1\relax
      \edef\setshift{%
        \noexpand\pgfkeys{
          /pgfplots/cycle list shift=\the\c@pgf@countb,
          /pgfplots/cycle list set=
        }
      }%
      \setshift%
      \globaldefs=0\relax
    \fi
    \pgfkeysgetvalue{/pgfplots/cycle list shift}\pgfplots@listindex@shift
    \ifx\pgfplots@listindex@shift\pgfutil@empty
    \else
      \c@pgf@counta=\pgfplots@listindex\relax
      \advance\c@pgf@counta by\pgfplots@listindex@shift\relax
      \ifnum\c@pgf@counta<0
        \c@pgf@counta=-\c@pgf@counta
      \fi
      \edef\pgfplots@listindex{\the\c@pgf@counta}%
    \fi
    \ifnum\pgfplots@cycle@dim>0
      \c@pgf@counta=\pgfplots@cycle@dim\relax
      \c@pgf@countb=\pgfplots@listindex\relax
      \advance\c@pgf@counta by-1
      \pgfplotsloop{%
        \ifnum\c@pgf@counta<0
          \pgfplotsloopcontinuefalse
        \else
          \pgfplotsloopcontinuetrue
        \fi
      }{%
        \pgfkeysgetvalue{/pgfplots/cycle multi list/@N\the\c@pgf@counta}\pgfplots@cycle@N
        \pgfplotsmathmodint{\c@pgf@countb}{\pgfplots@cycle@N}%
        \divide\c@pgf@countb by \pgfplots@cycle@N\relax
        \expandafter\pgfplots@getautoplotspec@
        \csname pgfp@cyclist@/pgfplots/cycle multi list/@list\the\c@pgf@counta @\endcsname
        {\pgfplots@cycle@N}%
        {\pgfmathresult}%
        \t@pgfplots@toka=\expandafter{#1,}%
        \t@pgfplots@tokb=\expandafter{\pgfplotsretval}%
        \edef#1{\the\t@pgfplots@toka\the\t@pgfplots@tokb}%
        \advance\c@pgf@counta by-1
      }%
    \else
      \pgfplotslistsize\autoplotspeclist\to\c@pgf@countd

      \pgfplots@getautoplotspec@{\autoplotspeclist}{\c@pgf@countd}{\pgfplots@listindex}%
      \let#1=\pgfplotsretval
    \fi
    \pgfmath@smuggleone#1%
    \endgroup
  }

\pgfplotsset{
  cycle list set/.initial=
}
\makeatother

\usepackage{lipsum}
\usepackage{amsfonts}
\usepackage{graphicx}
\usepackage{epstopdf}
\usepackage{amsmath,amssymb}
\ifpdf
  \DeclareGraphicsExtensions{.eps,.pdf,.png,.jpg}
\else
  \DeclareGraphicsExtensions{.eps}
\fi

\headers{Towards Matrix-Free Patch Smoothers for the Stokes Problem}{M. Wichrowski}

\title{Towards Matrix-Free Patch Smoothers for the Stokes Problem: Evaluating Local p-Multigrid Solvers}
\author{Michał Wichrowski}

\usepackage{amsopn}
\DeclareMathOperator{\diag}{diag}

\makeatletter
\newcommand*{\addFileDependency}[1]{
  \typeout{(#1)}
  \@addtofilelist{#1}
  \IfFileExists{#1}{}{\typeout{No file #1.}}
}
\makeatother

\providecommand{\keywordsname}{Keywords}
\newcommand{\Keywords}[1]{\par\noindent\textbf{\keywordsname:} #1}
\newcommand{\AMS}[1]{\par\noindent\textbf{AMS subject classifications:} #1}

\begin{document}
\author{Michał  Wichrowski$ ^0$}
\footnotetext{{Interdisciplinary Center for Scientific Computing, Heidelberg University, Germany. \texttt{mt.wichrowsk@uw.edu.pl}}}
\date{}

\maketitle

\begin{abstract}
    Vertex-patch smoothers offer an effective strategy for achieving robust geometric multigrid convergence for the Stokes
equations, particularly in the context of high-order finite elements. However, their practical efficiency is often
limited by the computational cost of solving the local saddle-point problems, especially when explicit matrix
factorizations are not feasible. We explore a fully iterative, matrix-free-compatible approach to the local patch solve
using $p$-multigrid techniques. We evaluate different local solver configurations: Braess-Sarazin and block-triangular
preconditioners. Our numerical experiments suggest that the Braess-Sarazin approach is particularly resilient. We find
that a single iteration of the local solver yields global convergence rates comparable to those obtained with exact
local solvers, even on distorted meshes and in the presence of large viscosity jumps.
\end{abstract}


\AMS{65Y10, 65Y20, 65N55, 65N30}

\section{Introduction}
\label{sec:intro}

Geometric multigrid methods establish the gold standard for elliptic partial differential equations, delivering optimal
$\mathcal{O}(N)$ complexity for symmetric positive definite problems~\cite{Hackbusch85, Bramble93, thomas2003textbook}.
However, extending this efficiency to the Stokes equations presents a fundamental structural challenge: the smoother
must not only reduce high-frequency errors in the velocity but also effectively damp error components within the
divergence-free subspace. Standard point-wise or cell-wise smoothers often fail in this regard because the support of
discrete divergence-free basis functions typically extends beyond a single element, preventing the local resolution of
the mass conservation constraint. Vertex-patch smoothers overcome this by solving coupled problems on overlapping
subdomains that fully enclose these supports, thereby recovering the stability of the splitting~\cite{hong2016robust,
  KanschatMao15}. The potential efficiency of this approach is close to direct solving the system; for instance, recent
work on GPU-accelerated solvers utilizing Raviart-Thomas elements~\cite{RaviartThomas, RaviartThomas2} demonstrated
global convergence in as few as 2 to 3 iterations~\cite{cui2025Stokes}. Yet, realizing this efficiency for general
high-order discretizations requires overcoming the significant computational bottleneck of solving these local patch
problems.

While high-order methods offer superior approximation properties and favorable dispersion error
characteristics~\cite{deville2002high}, they induce linear systems that are expensive to solve using traditional
matrix-based approaches. The memory bandwidth bottleneck associated with assembling and storing large, dense element
matrices has driven a paradigm shift toward matrix-free implementations~\cite{Kronbichler2012, Kronbichler2017a}. By
leveraging sum-factorization techniques on tensor-product cells, matrix-free methods reduce the operator evaluation
cost to $\mathcal{O}(p^{d+1})$ and minimize data movement, which is critical on modern hardware architectures
~\cite{kronbichler2019multigrid, munch2023cache}. However, this matrix-free constraint severely limits the choice of
available preconditioners and smoothers, as access to individual matrix entries is no longer possible.

To address the need for robust smoothing in high-order, matrix-free regimes, overlapping domain decomposition
techniques: specifically vertex-patch (or vertex-star) smoothers, have emerged as a powerful solution. Originally
analyzed in the context of domain decomposition~\cite{pavarino1993additive, ArnoldFalkWinther00}, these methods
decompose the global domain into small, overlapping subdomains centered around mesh vertices. By solving local problems
on these patches, they capture the local high-frequency spectrum of the operator, providing robust smoothing that is
independent of the mesh size $h$ and often the polynomial degree $p$~\cite{JanssenKanschat11, KanschatMao15}. This
approach has recently gained traction for complex problems, including non-matching grid
methods~\cite{wichrowski2025geometric, cui2025multigrid, bergbauer2025high} and $H(\text{div})$-conforming
discretizations~\cite{hong2016robust}. Moreover, the patch-based structure naturally promotes excellent data locality
with high arithmetic intensity, making it highly suitable for parallelization on both
CPUs~\cite{wichrowski2025smoothers,wichrowski2025pMGimplementaion} and GPUs~\cite{cui2025implementation}.

Despite these advantages, the adoption of patch smoothers is traditionally hindered by a single, critical bottleneck:
the computational cost of solving the local patch problems. For scalar elliptic problems, such as the Poisson equation,
highly efficient local solvers have been developed. These range from fast diagonalization methods (FDM) exploiting
tensor-product separability on structured grids~\cite{WitteArndtKanschat21, witte2025tensor} to specialized direct
solvers that utilize static condensation and Cholesky factorization~\cite{brubeck2021scalable}. However, for the Stokes
equations, there is no such \emph{silver bullet}. The coupling between velocity and pressure, combined with the Schur
complement, breaks the operator's separability even on perfectly Cartesian grids. Consequently, fast diagonalization
cannot be directly applied to the full system.

The literature offers several strategies to tackle the global Stokes system, which influence the design of local
solvers. A prominent class is block preconditioning~\cite{krzyzanowski2001block, krzyzanowski2011block}, where the
system is decoupled into velocity and pressure sub-problems. These methods typically employ approximations of the
pressure Schur complement, such as the pressure mass matrix or diffusion operators~\cite{olshanskii2006uniform}. While
effective as global preconditioners, their application within a local patch smoother has not been studied. The local
patch problems inherit the indefinite nature of the global system, and constructing robust block preconditioners that
respect the local boundary conditions of the patch -- especially on distorted meshes -- remains an open challenge.

An alternative approach, specifically designed for multigrid smoothing, is the class of coupled \emph{Vanka} or
box-relaxation smoothers~\cite{vanka1986block}. These methods solve small, local saddle-point systems associated with a
cell or a cluster of degrees of freedom. Vanka smoothers have been extensively studied and applied to various
discretizations~\cite{wobker2009numerical, farrell2021local}, including variable viscosity
problems~\cite{borzacchiello2017box} and space-time formulations~\cite{margenberg2025hp, margenberg2025multigrid}.
While classical Vanka smoothers often rely on dense local inverses or direct solvers, recent efforts have explored
inexact or compressed variants to reduce computational load~\cite{harper2023compression}. However, ensuring the
$p$-robustness of these coupled smoothers without incurring excessive cost remains a difficulty, particularly for
matrix-free high-order methods where forming the local dense matrix comes at a prohibitive cost.

The robustness of the local solver becomes paramount when targeting complex physical applications characterized by
large jumps in material coefficients. Two prime examples are Geodynamics, specifically mantle
convection~\cite{rudi2017weighted, burstedde2011p4est}, and Fluid-Structure Interaction
(FSI)~\cite{wichrowski2023exploiting, voronin2025monolithic}. In FSI problems solved via monolithic schemes, the solid
subdomain behaves as a fluid with extremely high viscosity (pseudo-viscosity proportional to the inverse time step),
leading to viscosity contrasts of orders $10^6$ to $10^{12}$ across the interface. Similarly, mantle convection
involves viscosity variations of several orders of magnitude due to temperature and pressure dependence. In such
regimes, standard smoothers often stall. The local patch solver must therefore be robust not only to mesh distortion
but also to extreme coefficient heterogeneity, effectively capturing the local interface
dynamics~\cite{olshanskii2006uniform, may2015scalable}.

A promising candidate for handling these difficult regimes is the Braess-Sarazin smoother~\cite{braess1997efficient}.
As elegantly described in~\cite{chen2015multigrid}, this method can be interpreted as a constrained smoother that
effectively handles the saddle-point structure by maintaining the velocity-pressure coupling. While it has not been the
standard choice in the broader community, recent investigations have demonstrated its specific efficacy for matrix-free
implementations~\cite{wichrowski2022matrix, jodlbauer2024matrix}. Furthermore, it has been shown to perform reasonably
well for the extreme viscosity contrasts inherent in Fluid-Structure Interaction
problems~\cite{wichrowski2023exploiting}. However, utilizing Braess-Sarazin as a \textit{local} solver inside a
vertex-patch smoother implies a nested complexity: one must solve a saddle-point problem inside a patch, which itself
is part of a global multigrid hierarchy.

In this work, we embrace this recursive spirit to overcome the local solver bottleneck. Instead of relying on expensive
exact direct solvers or restrictive separable approximations, we propose constructing the local patch solver using a
$p$-multigrid method (coarsening in polynomial degree)~\cite{wichrowski2025local}. This approach applies the logic of
geometric multigrid one level deeper: the local patch problem is solved iteratively by traversing a hierarchy of
polynomial orders. This strategy preserves the matrix-free, sum-factorized structure of the operator at all levels,
maintaining optimal memory complexity and allowing for the reuse of the global matrix-free
infrastructure~\cite{dealii2019design, dealII97}.

We investigate and compare three distinct iterative strategies for the local Stokes patch problem, all underpinned by
this efficient $p$-multigrid solver for the velocity block:
\begin{enumerate}
  \item \textbf{Local Braess-Sarazin:} We adapt the Braess-Sarazin smoother to act as the local solver on the patch. Crucially, we employ $p$-multigrid to approximate the velocity inverse required within the Braess-Sarazin iteration. This effectively nests a block-smoother inside a patch-smoother inside a global multigrid cycle.
  \item \textbf{Local Block GMRES:} We utilize a Krylov subspace method preconditioned with a block-triangular operator on the patch, again using $p$-multigrid for the velocity block.
\end{enumerate}
Additionally, for Cartesian patches, we explore a block-by-block elimination approach, solving for the pressure via the Schur complement, similar to recent work on GPUs~\cite{cui2025Stokes}, but substituting exact inverses with $p$-multigrid approximations.  Noteworthy, a high-order system can also be preconditioned by a refined mesh of degree-one elements~\cite{pazner2020efficient}, but this approach is outside the scope of our current matrix-free framework.

Our numerical results demonstrate that replacing exact local solvers with iterative $p$-multigrid approximations does
not compromise the robustness of the global multigrid method. Specifically, the Braess-Sarazin variant proves to be
exceptionally robust, handling distorted meshes and large viscosity jumps with iteration counts comparable to using
exact local inverses.

The objective of this work is to establish the mathematical robustness and convergence properties of the local
p-multigrid smoother. While a high-performance matrix-free implementation is the ultimate goal, this paper focuses on
the algorithmic foundation and demonstrates that the 'inexact' nature of the local solver does not degrade the global
convergence rate --- a prerequisite for any future optimized implementation.

The paper is organized as follows. Section~\ref{sec:multigrid} provides the mathematical background for the Stokes
problem and vertex patch smoothers. Section~\ref{sec:local-stokes-solvers} describes the local solvers in detail,
including the matrix-free $p$-multigrid framework. Section~\ref{sec:results} presents numerical results comparing the
performance of the different local solvers. For the reader interested in key performance metrics, we highlight the
comparison of local solver variants in Table~\ref{tab:single_patch_cartesian}, and the global multigrid convergence
results demonstrating robustness against mesh distortion and high-contrast viscosity jumps in
Table~\ref{tab:gmres_iterations_jacobi_2D} and Table~\ref{tab:gmres_iterations_mu1e6_2D}, respectively. Finally,
Section~\ref{sec:conclusion} concludes the paper.
\section{Mathematical and Algorithmic Background}
\label{sec:multigrid}
Let us first briefly describe the mathematical setup of the multigrid method, first
assume a hierarchy of meshes
\begin{gather}
  \mesh_0 \sqsubset \mesh_1 \sqsubset \dots \sqsubset \mesh_L,
\end{gather}
subdividing a domain in $\R^d$, where the symbol
``$\sqsubset$'' indicates nestedness, that is, every cell of mesh
$\mesh_{\ell+1}$ is obtained from a cell of mesh $\mesh_{\ell}$ by refinement. We
note that topological nestedness is sufficient from the algorithmic
point of view, such that domains with curved boundaries can be
covered approximately.

For our numerical experiments, we consider the Stokes equations on a domain $\Omega \subset \mathbb{R}^d$. The weak
formulation reads: find velocity $u \in V = [H^1_0(\Omega)]^d$ and pressure $p \in Q = L^2_0(\Omega)$ such that
\begin{gather}
  \label{eq:bilinear-form}
  a((u,p),(v,q)) = (\nabla u, \nabla v) - (\nabla \cdot v, p) + (\nabla \cdot u, q) = (f, v) \quad \forall (v,q) \in V \times Q.
\end{gather}

To discretize the problem, we associate each mesh $\mesh_\ell$ with a pair of finite element spaces $(\mathbb{V}_\ell,
  \mathbb{Q}_\ell) \subset \mathbb{V} \times \mathbb{Q}$. Specifically, we use the $\mathbb{Q}_p-\mathbb{P}_{p-1}$ pair,
which consists of continuous polynomials of degree $p$ for the velocity and discontinuous polynomials of degree $p-1$
for the pressure. We associate each space with its coefficient vector $\mathbb R^{\operatorname{dim}V_\ell}$ and
$\mathbb R^{\operatorname{dim}Q_\ell}$ and do not distinguish between a finite element function and its coefficient
vector by notation. Between these spaces, we introduce transfer operators
\begin{gather}
  \begin{array}{rlcl}
    \resth\ell\colon & V_{\ell+1} & \to & V_\ell,     \\
    \prolh\ell\colon & V_{\ell}   & \to & V_{\ell+1}.
  \end{array}
\end{gather}
As usual, $\prolh\ell$ is chosen as the embedding operator and $\resth\ell$ as its
$\ell_2$-adjoint. The subscript $h$ indicates h-multigrid transfer between meshes of different resolution.

Note that while the finite element spaces are nested ($V_\ell \subset V_{\ell+1}$ and $Q_\ell \subset Q_{\ell+1}$), the
resulting discretely divergence-free subspaces are generally non-nested. Indeed, the Stokes problem can be viewed as a
symmetric positive definite problem when restricted to the space of divergence-free functions~\cite{chen2015multigrid}.
However, because our element pair does not enforce the divergence-free constraint pointwise (only a discrete
divergence-free condition), the corresponding subspaces on different levels are non-nested. This is because a function
being discretely divergence-free on level $\ell$ does not automatically satisfy the more restrictive weak divergence
condition on level $\ell+1$. True nestedness of the divergence-free subspaces is only possible when the discretization
achieves an exact pointwise divergence-free condition, as is the case for certain $H(\text{div})$-conforming
discretizations~\cite{RaviartThomas,RaviartThomas2, ArnoldFalkWinther00}.

The discrete problem on level $\ell$ reads
\begin{gather}
  \label{eq:matrix}
  \mathcal{A}_\ell x_\ell = \begin{pmatrix} A_\ell & B_\ell^T \\ B_\ell & 0 \end{pmatrix} \begin{pmatrix} u_\ell \\ p_\ell \end{pmatrix} = \begin{pmatrix} f_\ell \\ 0 \end{pmatrix},
\end{gather}
where $\mathcal{A}_\ell$ is the block system matrix obtained by discretization.

Multigrid methods are common solution methods for the discrete linear system \eqref{eq:matrix}. The efficiency of the
overall method, however, depends on the choice of the so-called smoother. For saddle point problems like Stokes,
constructing effective smoothers is challenging. In this work, we focus on a vertex patch smoother, which has been
shown to yield robust convergence. We describe this smoother in detail in the next subsection.

\subsection{Vertex patch smoothers}
The vertex patch smoother is an overlapping subspace correction
method~\cite{ArnoldFalkWinther00,JanssenKanschat11,KanschatMao15,WitteArndtKanschat21,brubeck2021scalable,wichrowski2025smoothers}.
The domain is decomposed into a collection of overlapping subdomains, or \emph{patches} $\Omega_j$, each formed by the
cells surrounding an interior vertex. The smoother iteratively improves the solution by solving local problems on these
patches. We employ a multiplicative (Gauss-Seidel-like) variant, where patches are processed sequentially.

The update for a single patch $j$ consists of four main steps. To define these, we introduce restriction operators:
$\overline{\Pi}_j$ extracts all degrees of freedom (DoFs) on a patch $\Omega_j$ (including its boundary) from a global
vector, while $\Pi_j$ extracts only the DoFs interior to the patch.
\begin{enumerate}
  \item \textbf{Gather.} Collect the current solution values on the patch into a local vector: $\overline{u}_j =
          \overline{\Pi}_j u$.
  \item \textbf{Evaluate.} Compute the local residual $r_j$ for the interior DoFs. This uses the local operator
        $\overline{A}_j$ (acting on all patch DoFs) and the global right-hand side $b$: $r_j = \Pi_j b - \Pi_j
          \overline{A}_j \overline{u}_j$.
  \item \textbf{Local solve.} Find a correction $d_j$ by approximately solving the local system on the patch
        interior: $d_j \approx \tilde A_j^{-1} r_j$, where $\tilde A_j^{-1}$ is an inexpensive approximate solver.
  \item \textbf{Scatter.} Add the local correction $d_j$ back to the global solution vector using the transpose of the
        interior restriction operator: $u \pluseq \Pi_j^T d_j$.
\end{enumerate}
These steps are summarized in Algorithm~\ref{alg:Loop-sequential}. This implementation follows the cell-oriented loop
structure proposed in~\cite{wichrowski2025smoothers}, which is well-suited for matrix-free operator evaluation.

\begin{algorithm}[tp]
  \begin{algorithmic}
    \State \Function{\LUpNew}{j,u}
    \State    $ u_j \gets \overline\Pi_j u$ \Comment*{Gather}
    \State    $ r_{j} \gets \Pi_j b - \Pi_j\overline A_j u_j$ \Comment*{Evaluate}
    \State    $ d_j \gets \tilde{A}_j ^{-1} r_j$ \Comment*{Solve}
    \State \Return{ $\Pi_j^T d_j$ }\Comment*{Scatter}
    \EndFunction
    \State \For{$j=1,\ldots,N_\text{patches}$}	\Comment{Main smoother sweep}
    \State $u \pluseq$ \Call{\LUpNew}{$j$,$u$}
    \EndFor
  \end{algorithmic}
  \caption{Application of a patch smoother where local residuals are computed on-the-fly and combined with a local
    solver.}
  \label{alg:Loop-sequential}
\end{algorithm}

In this work, we focus on the local solve (step 3), which is critical for the overall performance. This step, denoted
as the application of an approximate inverse $\tilde{A}_j^{-1}$ in Algorithm~\ref{alg:Loop-sequential}, is performed by
applying a few iterations of a preconditioned Richardson method, as outlined in Algorithm~\ref{alg:local-richardson}.
Notice that the first iteration is performed outside the loop to avoid an unnecessary matrix-vector product, since the
initial guess for the correction is zero.

\begin{algorithm}[htbp]
  \caption{Local preconditioned Richardson iteration for the patch problem. This function implements the action of
  $\tilde{A}_j^{-1}$ from Algorithm~\ref{alg:Loop-sequential}.}
  \label{alg:local-richardson}
  \begin{algorithmic}[1]
    \Function{LocalSolve}{$\residual_j$}
    \Comment{The input $r_j$ is the residual for the patch problem.}
    \State $d_j \gets \Call{p-V-cycle}{\residual_j}$ \Comment{First iteration with zero initial guess.}
    \For{$k=2,\ldots,N_\text{iter}$} \Comment{Apply remaining steps.}
    \State $\res \gets \residual_j - A_j d_j$ \Comment{Compute residual for the correction.}
    \State $c_j \gets \Call{p-V-cycle}{\res}$ \Comment{Precondition with one p-MG cycle.}
    \State $d_j \pluseq c_j$ \Comment{Update correction.}
    \EndFor
    \State \Return{$d_j$}
    \EndFunction
  \end{algorithmic}
\end{algorithm}

\begin{figure}[tp]
  \centering
\def\svgwidth{\columnwidth}
\begin{tikzpicture}[scale=0.9, transform shape, font=\normalsize,
        cell/.style={draw=black, thick, minimum size=1.2cm},
        interior/.style={fill=blue!20},
        global/.style={fill=red!15},
        panel-label/.style={font=\bfseries\normalsize, align=center},
        arrow/.style={->, thick, shorten >=2mm, shorten <=2mm},
        dof/.style={circle, fill=black, inner sep=0.8pt}
    ]

    \newcommand{\drawPatchInteriorDoFs}[2]{
        \foreach \i in {1,...,5} {
                \foreach \j in {1,...,5} {
                        \pgfmathsetmacro{\xi}{\i/6}
                        \pgfmathsetmacro{\eta}{\j/6}

                        \coordinate (NW) at (#1.north west);
                        \coordinate (SE) at (#2.south east);

                        \coordinate (pt_x) at ($(NW)!\xi!(SE)$);	 
                        \coordinate (pt_y) at ($(NW)!\eta!(SE)$);  
                        \coordinate (pt) at (pt_x |- pt_y);
                        \node[dof, minimum size=4pt, inner sep=0pt, fill=blue] at (pt) {};
                    }
            }
    }

    \newcommand{\drawCubicDoFs}[1]{
        \foreach \xi in {0, 0.33, 0.66, 1} {
                \foreach \eta in {0, 0.33, 0.66, 1} {
                        \coordinate (bottom) at ($(#1.south west)!\xi!(#1.south east)$);
                        \coordinate (top) at ($(#1.north west)!\xi!(#1.north east)$);
                        \node[dof, minimum size=3pt, inner sep=0pt,  fill=red] at ($(bottom)!\eta!(top)$) {};
                    }
            }
    }

    \newcommand{\drawBoundaryDoFs}[2][dof]{
        \foreach \xi in {0, 0.33, 0.66, 1} {
                \coordinate (pt) at ($(#2.south west)!\xi!(#2.south east)$);
                \node[#1] at (pt) {};
            }
        \foreach \xi in {0, 0.33, 0.66, 1} {
                \coordinate (pt) at ($(#2.north west)!\xi!(#2.north east)$);
                \node[#1] at (pt) {};
            }
        \foreach \eta in {0.33, 0.66} {
                \coordinate (pt) at ($(#2.south west)!\eta!(#2.north west)$);
                \node[#1] at (pt) {};
            }
        \foreach \eta in {0.33, 0.66} {
                \coordinate (pt) at ($(#2.south east)!\eta!(#2.north east)$);
                \node[#1] at (pt) {};
            }
    }
    \newcommand{\drawInteriorFrame}[2]{%
        \draw[rounded corners=3pt, thick, draw=blue!60!black, fill=gray!80, fill opacity=0.1, densely dashed]
        ([shift={(4pt,-4pt)}]#1.north west)
        rectangle
        ([shift={(-4pt,4pt)}]#2.south east);%
    }

    \node[cell] (p1c1) at (-0.2,0) {};
    \node[cell] (p1c2) at (1,0) {};
    \node[cell] (p1c3) at (-0.2,-1.2) {};
    \node[cell] (p1c4) at (1,-1.2) {};
    \drawInteriorFrame{p1c1}{p1c4}

    \drawCubicDoFs{p1c1}
    \drawCubicDoFs{p1c2}
    \drawCubicDoFs{p1c3}
    \drawCubicDoFs{p1c4}
    \drawPatchInteriorDoFs{p1c1}{p1c4}

    \node[cell] (p2c1) at (4.5,0) {}; 
    \node[cell] (p2c2) at (5.7,0) {}; 
    \node[cell] (p2c3) at (4.5,-1.2) {}; 
    \node[cell] (p2c4) at (5.7,-1.2) {}; 
    \drawInteriorFrame{p2c1}{p2c4}
    \drawPatchInteriorDoFs{p2c1}{p2c4}

    \node[cell] (p3c1) at (9,0) {}; 
    \node[cell] (p3c2) at (10.2,0) {}; 
    \node[cell] (p3c3) at (9,-1.2) {}; 
    \node[cell] (p3c4) at (10.2,-1.2) {}; 
    \drawInteriorFrame{p3c1}{p3c4}
    \drawPatchInteriorDoFs{p3c1}{p3c4}

    \draw[arrow] ($(p1c3.west) - (2cm, 0)$) --	node[above, midway, ] {$\overline\Pi_j u, \;\;	\Pi_j b $} node[below,
        midway, ] {Gather}  (p1c3.west);
    \draw[arrow] (p1c4.east) -- node[above, midway, ] {$b_j - \Pi_j\overline A_j \overline u_j$} node[below, midway, ]
    {Evaluate} (p2c3.west);
    \draw[arrow] (p2c4.east) -- node[above, midway, ]{$\tilde A_j^{-1} r_j$}
    node[below, midway, ] {Local Solve}   (p3c3.west);
    \draw[arrow] (p3c4.east) -- node[above, midway, ]{$\Pi_j^T d_j$} node[below, midway, ] {Scatter} ($(p3c4.east) +
    (2cm, 0)$);

\end{tikzpicture}
  \caption{Workflow in a generic smoother application of a patch smoother for a single patch $j$. The three panels
    represent the main steps:
    gathering data from the global solution and right-hand side, computing the local residual, and applying the local
    correction and scattering it back to the global solution. Blue dots indicate interior DoFs of the patch, while red
    dots indicate all DoFs (including boundary) associated with each cell. The dashed rectangle highlights the patch
    interior. Figure taken from~\cite{wichrowski2025pMGimplementaion}.}
  \label{fig:patch-smoother-steps}
\end{figure}
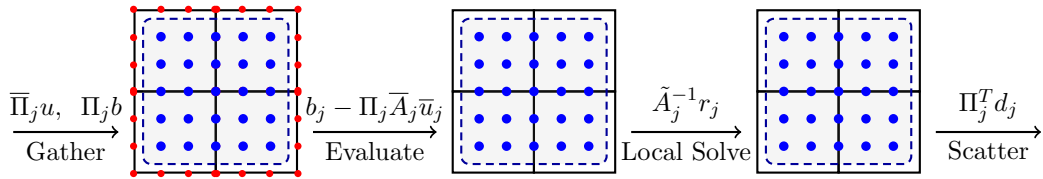

\subsection{Local Solvers for Stokes Patch Problems}
\label{sec:local-stokes-solvers}
The efficiency of the patch smoother hinges on the local solve step. For the Stokes equations, this involves solving a
saddle point system on each patch. We consider three iterative approaches for this task: the Braess-Sarazin smoother, a
block GMRES solver, and a block-by-block approach using the Schur complement.

All these methods rely on a local multigrid solver. Since a patch contains too few cells for effective geometric
coarsening, we construct a hierarchy of levels by reducing the polynomial degree $p$. This p-multigrid approach,
detailed in~\cite{wichrowski2025local}, provides an effective and computationally efficient way to approximate the
action of the inverse. It is computationally inexpensive to apply while keeping memory requirements low by relying on
matrix-free data structures and storing only compact per-patch data.

We distinguish between two main ways to employ this local p-multigrid solver:
\begin{enumerate}
  \item \textbf{Multigrid on the Stokes System (Braess-Sarazin):} The p-multigrid V-cycle is applied to the full
        Stokes system $M_j$. The smoother on each level is the Braess-Sarazin iteration.
  \item \textbf{Multigrid on the Velocity Block (Block Solvers):} The p-multigrid V-cycle is applied only to the
        velocity block $A_j$ to approximate its inverse. This approximate inverse is then used within a block
        preconditioner for the full system.
\end{enumerate}

\subsubsection{Local p-Multigrid Framework}
To solve a local patch problem $M_j d_j = r_j$, where $M_j$ is either the full Stokes system or the velocity block
$A_j$, we employ a p-multigrid method. As with any multigrid method, the p-multigrid solver is built upon a few key
components: a hierarchy of levels, transfer operators between these levels, a smoother $S_{j,p}$ (with relaxation
parameter $\omega$ and $m$ steps) and operator $M_{j,p}$ on each level, and a coarse-grid solver.

In the p-multigrid context, the levels correspond to a hierarchy of polynomial degrees. We choose a geometric
progression, such as $p_l=1, 3, 7, \dots$, up to the target degree. If the target degree does not fit this sequence,
the hierarchy is adjusted accordingly; for instance, for a target degree of $p=4$, the sequence would be $1, 3, 4$. The
transfer operators are the natural embedding for prolongation (from degree $p_{l-1}$ to $p_l$) and its adjoint for
restriction. The explicit Kronecker product structure of these operators and the closed-form definitions of the level
operators are detailed in Appendix~\ref{sec:appendix-local-inverse}.

The complete p-multigrid V-cycle is outlined in Algorithm~\ref{alg:p-v-cycle} in which we also include the option for
skipping the post-smoothing step. The algorithm recursively descends to the coarsest level, where an exact solve is
performed.

\begin{algorithm}[htbp]
  \caption{The p-multigrid V-cycle for a local patch problem. The post-smoothing step (lines 12-13) is omitted in the
    \emph{half V-cycle} optimization.}
  \label{alg:p-v-cycle}
  \begin{algorithmic}[1]
    \Function{V-cycle}{$l, \residual_p$}
    \If{$l = 1$}
    \State $d_1 \gets M_{j,1}^{-1} \residual_1$ \Comment{Coarse-level solve}
    \State \Return $d_1$
    \EndIf
    \State $d_p \gets \omega S_{j,p}^{-1} \residual_p$ \Comment{Pre-smoothing}
    \State $\res_p \gets \residual_p - M_{j,p} d_p$ \Comment{Compute residual}
    \State $\res_{p-1} \gets \restp{p,p-1} \res_p$ \Comment{Restrict residual}
    \State $e_{p-1} \gets \Call{V-cycle}{p-1, \res_{p-1}}$ \Comment{Recursive call}
    \State $e_p \gets \prolp{p-1,p} e_{p-1}$ \Comment{Prolongate correction}
    \State $d_p \gets d_p + e_p$ \Comment{Apply correction}
    \If{full V-cycle}
    \State $\res_p \gets \residual_p - M_{j,p} d_p$ \Comment{Compute residual for post-smoothing }
    \State $d_p \gets d_p + \omega S_{j,p}^{-1} \res_p$ \Comment{Post-smoothing }
    \EndIf
    \State \Return $d_p$
    \EndFunction
  \end{algorithmic}
\end{algorithm}

\begin{figure}[tp]
  \input{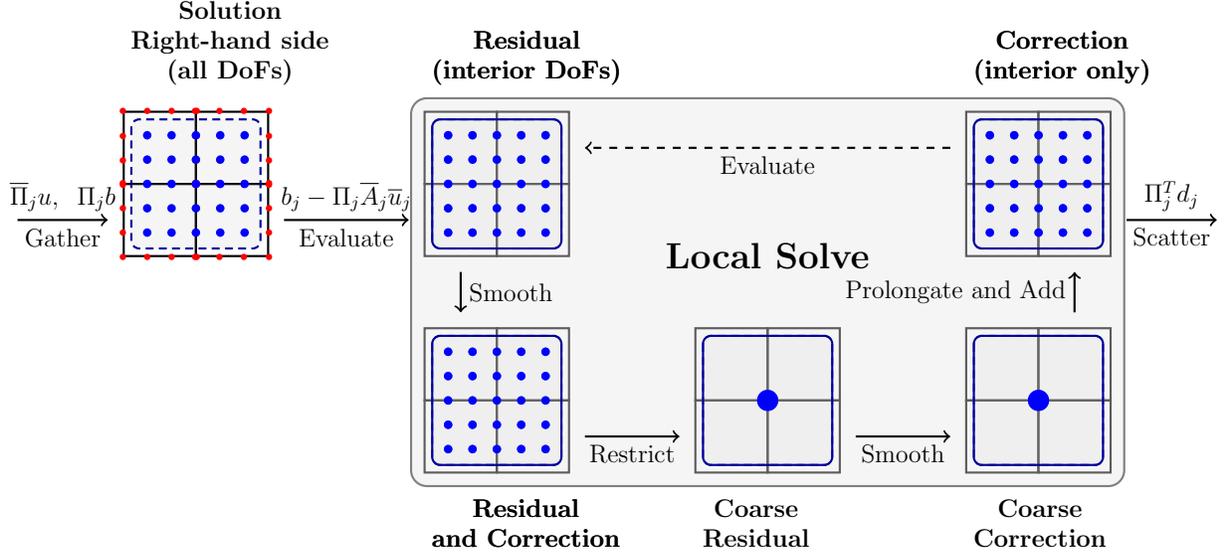}
  \caption{Patch-smoother data flow for a local p-multigrid solve (illustration for $p=3$ on the fine local level). The
    coarse level shown is $p=1$ and contains a single interior node. The
    diagram
    shows gather, evaluate (local residual), pre-smoothing, restriction to the coarse level, an exact coarse
    correction, prolongation-and-add, and scatter-add. A dashed arrow indicates proceeding to the next local
    multigrid iteration. For brevity we illustrate only pre-smoothing; post-smoothing is omitted. Figure taken from~\cite{wichrowski2025pMGimplementaion}.}
  \label{fig:patch-flow}
\end{figure}

\subsubsection{Multigrid on the Stokes System}
\label{sec:bs-smoother}
In this approach, we apply multigrid to the entire Stokes problem on the patch. As with any geometric or p-multigrid method, the efficiency of the solver relies on three ingredients: the transfers and level operators, the coarse grid solver, and the smoother.

The definition of the level systems and transfer operators is standard. The level operators correspond to the
discretization of the Stokes problem using the finite element space of polynomial degree $p$ associated with the
current level. The transfer operators are the natural embeddings between the nested velocity and pressure spaces and
their adjoints.

The coarse grid solver, however, requires special consideration. The hierarchy of levels bottoms out at $p=1$. While
the velocity space is defined on the mesh of the patch, we define the pressure space on a coarsened geometry; treating
the entire patch as a single macro-cell (e.g., a single square instead of the four cells that typically form a 2D
patch). This choice results in a pressure space consisting of only a single constant degree of freedom over the whole
patch. Since we solve local problems with homogeneous Dirichlet boundary conditions for the velocity, the pressure is
only defined up to a constant. To ensure a unique solution, we fix the average pressure on the patch to zero. However,
enforcing a zero mean on a space of constants constrains the single available degree of freedom to be zero.
Consequently, the effective pressure space on the coarsest level is empty, and the coarse grid solve reduces to a
symmetric positive definite problem for the velocity variables only.

For the smoother $S_{j,p}$ on the levels $p > 1$, we employ the Braess-Sarazin smoother~\cite{braess1997efficient}. Use
of this smoother is motivated by its robustness with respect to the stationary or non-stationary nature of the flow.
Each smoothing step involves solving the saddle point system
\begin{gather}
  \begin{pmatrix} \tilde{A}_j & B_j^T \\ B_j & 0 \end{pmatrix} \begin{pmatrix} \delta u_j \\ \delta p_j \end{pmatrix} = \begin{pmatrix} r_{u,j} \\ r_{p,j} \end{pmatrix},
\end{gather}
where $\tilde{A}_j$ is a symmetric positive definite approximation of the velocity block $A_j$.
This system is solved using block-wise elimination, which proceeds in three stages:
\begin{enumerate}
  \item Compute an intermediate velocity correction: $\delta u_j^{(1)} = \tilde{A}_j^{-1} r_{u,j}$.
  \item Solve the Schur complement system for the pressure correction $\delta p_j$:
        \begin{gather}
          \label{eq:bs-schur}
          S_j \delta p_j = B_j \delta u_j^{(1)} - r_{p,j}, \quad \text{with } S_j = B_j \tilde{A}_j^{-1} B_j^T.
        \end{gather}
  \item Perform the final velocity update: $\delta u_j = \delta u_j^{(1)} - \tilde{A}_j^{-1} B_j^T \delta p_j$.
\end{enumerate}
The resulting smoothing operator is denoted by $\mathcal{S}_k^{-1}$ and its closed-form block representation is
presented in Appendix~\ref{sec:appendix-local-inverse}.

The effectiveness of the Braess-Sarazin smoother relies on one of two conditions\cite{zulehner2000class}: either
$\tilde{A}_j$ is a high-quality preconditioner for $A_j$, or the Schur complement system~\eqref{eq:bs-schur} is solved
accurately. In this work, we follow the latter approach. The classical method proposed by
Zulehner~\cite{zulehner2000class} ensures robustness by solving the Schur complement system with a multigrid method.
However, it has been shown in~\cite{wichrowski2022matrix} that even a single iteration of a standard iterative solver
like the Conjugate Gradient (CG) method can be sufficient to maintain robust convergence, which significantly reduces
the computational overhead of the local solve.

We evaluate two options for the inner solver used to solve the Schur complement system \eqref{eq:bs-schur}: CG and
Richardson iteration. We denote the number of inner iterations by $n_S$. While CG is typically more efficient in terms
of iteration counts, it is a non-linear operator. Its use within the smoother requires the outer Krylov solver to be
flexible (e.g., FGMRES). In contrast, Richardson iteration is a linear operator and preserves the consistency of the
preconditioner. Furthermore, Richardson iteration offers a specific computational benefit: for a zero initial guess,
the first iteration does not require an evaluation of the Schur operator $S_j$. Since $S_j$ involves two applications
of $B$ and one of $\tilde{A}_j^{-1}$, skipping this evaluation for the first step provides a meaningful speedup. When
Richardson iteration is used, its relaxation parameter is autotuned based on an eigenvalue estimate of the
preconditioned operator to ensure stability and rapid convergence.

Both inner solvers are preconditioned with an approximation $\tilde{S}_j^{-1}$ of the inverse Schur complement, which
we typically choose as the (approximate) diagonal of $S_j$. If $\tilde{A}_j$ is chosen as the diagonal of $A_j$ (damped
Jacobi), the exact diagonal of $S_j$ can be computed. Even if $\tilde{A}_j$ is not diagonal, an effective
preconditioner can be constructed by approximating $\tilde{A}_j^{-1}$ with $(\text{diag}(A_j))^{-1}$ when computing the
Schur diagonal, as shown in~\cite{wichrowski2022matrix}. Although more advanced approximations for $\tilde{A}_j^{-1}$
are possible (e.g., using Chebyshev polynomials or scaling with eigenvalue estimates), we focus here on simple diagonal
approximations.

We note a special case for Cartesian grids where the Fast Diagonalization Method (FDM) allows for the exact inversion
of $A_j$ at low cost. In this scenario, one could technically use $\tilde{A}_j = A_j$ combined with a less accurate
Schur solver. This would resemble a block preconditioner approach but applied within the multilevel hierarchy. However,
since our focus is on general unstructured patches where FDM is not available, we do not pursue this direction further.

Additionally, the choice for the number of smoothing steps $m$ per level is motivated by theory. The analyses
in~\cite{braess1997efficient,zulehner2000class} are formulated for nested multigrid cycles with at least two smoothing
steps (pre- and post-smoothing) per level, i.e., $m=2$; these hypotheses underpin the theoretical robustness results.
In practice, however, one can try using a single smoothing step ($m=1$) to reduce cost. Therefore we explicitly test
both variants in our experiments: the two-step (classical) smoother that matches the theory, and the single-step
(economized) variant that omits post-smoothing. The single-step variant may weaken the strict theoretical guarantees
but can offer substantial runtime savings, and in our numerical tests it typically leads to only a modest increase in
outer iteration counts for the Braess--Sarazin approach.

An additional advantage of the Braess--Sarazin smoother is its robustness with respect to the specific form of the
equation, maintaining effectiveness for stationary Stokes, non-stationary variants, and even linearized operators like
the Oseen equations~\cite{wichrowski2023exploiting}. While convergence rates may vary depending on the flow regime, the
smoother's block-wise structure naturally balances mass and stiffness contributions, whereas some block preconditioners
can be more sensitive to the relative weighting of these terms~\cite{olshanskii2006uniform}.

\subsubsection{Block Solvers: Multigrid on the Velocity Block}
Alternatively, we can employ an iterative solver, such as GMRES with a block preconditioner, for the local Stokes
system. In this case, the p-multigrid method is used solely to approximate the inverse of the velocity block $A_j$. To
ensure symmetry and robustness, we use a symmetric block preconditioner $\mathcal{P}_j$ that mirrors the structure of
the Braess-Sarazin iteration:
\begin{gather}
  \mathcal{P}_j^{-1} = \begin{pmatrix} \tilde{A}_j & B_j^T \\ B_j & 0 \end{pmatrix}^{-1}
\end{gather}
where $\tilde{A}_j^{-1}$ is a high-quality approximation of the velocity inverse, typically obtained by one or more p-multigrid V-cycles (setting $M_j = A_j$ in Algorithm~\ref{alg:p-v-cycle}).
Similarly, $\tilde{S}_j^{-1}$ is a \emph{good} approximation of the inverse Schur complement, for instance, the action of a pressure mass matrix inverse or a simple iterative solve.

The advantage of GMRES is its ability to handle the indefinite nature of the saddle point system without requiring an
exact Schur complement solve. However, the downside is that the structure of the preconditioner is problem-depedent,
and with variable coefficients or non-stationary flows, its performance can degrade. In that case, a better
approxiamtion of the Schur complement can obtained by including a weighted inverse mass matrix and inverse laplacial
operator~\cite{olshanskii2006uniform}. Since laplacian inverses are not easily available in our patch setting, another
p-multigrid V-cycle could be used to approximate the laplacian inverse on the pressure space. However, this would
significantly increase the computational cost of the local solve, and we do not pursue this direction further.

\subsubsection{Blockwise Elimination via Schur Complement}
Finally, we consider a block-by-block iteration based on the Schur complement. This method follows the same three-stage
elimination procedure described for the Braess-Sarazin smoother in Section~\ref{sec:bs-smoother}, with the specific
choice $\tilde{A}_j = A_j$, i.e., using the exact inverse of the velocity block. This approach, which has also been
successfully utilized in~\cite{cui2025Stokes}, is primarily feasible on Cartesian grids where the Fast Diagonalization
Method (FDM)~\cite{Lynch1964} allows for the exact inversion of $A_j$ at a computational cost comparable to a few
matrix-vector products.

A significant challenge, even in this ideal Cartesian setting, remains the lack of a \emph{silver bullet} for the Schur
complement $S_j = B_j A_j^{-1} B_j^T$ comparable to FDM for the Poisson equation. While the Laplacian on a Cartesian
patch can be decomposed into a Kronecker product of 1D operators to enable an $O(N^{d+1/d})$ solve, the Stokes Schur
complement does not inherit this separable structure. The coupling between velocity components and the divergence
constraint breaks the separability required for FDM. Consequently, one must either rely on the explicit formation of
the (dense) Schur complement matrix which is only practical for low polynomial degrees, or use iterative solvers to
obtain the pressure. In~\cite{cui2025Stokes} an average of 15 CG iteration on Schur complement was reported. For
general unstructured or distorted patches where FDM is unavailable, the requirement of an exact velocity inverse would
make this block-by-block approach significantly more expensive than the standard Braess-Sarazin iteration, further
motivating our use of p-multigrid to approximate the velocity inverse.

\subsubsection{Summary of local solver variants}
The two iterative approaches Presented above (Braess--Sarazin and block solvers) share a common architectural framework
and are configured by the following tuning parameters:
\begin{itemize}
  \item the \textbf{number of smoothing steps} $m$ and the \textbf{relaxation parameter} $\omega$ used within the local
        p-multigrid levels (applied to either the full Stokes system or the velocity block);
  \item the \textbf{velocity block approximation} $\tilde{A}_j^{-1}$, ranging from a simple diagonal or damped Jacobi step to
        p-multigrid V-cycles or an exact inverse; in case if p-multigrid is used, then relaxation parameter $\omega$ and number
        of smoothing steps $m$ per level must be specified.
  \item the \textbf{Schur complement approximation} $\tilde{S}_j^{-1}$, typically chosen as a (preconditioned) diagonal or a
        pressure mass matrix inverse $\mathcal{M}_p^{-1}$.
\end{itemize}
For Braess--Sarazin  the following additional parameters are considered:
\begin{itemize}
  \item the \textbf{Schur solver}, either preconditioned CG or Richardson iteration;
  \item the \textbf{number of inner iterations} $n_S$ for the Schur complement solve.
\end{itemize}
Unless stated otherwise, we use $m=1$ smoothing step per level and $n_S=1$ inner Richardson iteration for the Schur complement in the Braess-Sarazin solver. We emphasize that when Richardson iteration is used for the Schur system, the relaxation parameter is autotuned
via an eigenvalue estimate of the preconditioned operator $\tilde{S}_j^{-1} S_j$.
In all experiments, we choose the damping parameter $\omega = 0.7$.
Interestingly,  the approximate local inverse obtained via a single p-multigrid V-cycle with one smoothing step per level can be written in closed form as a combination of Kronecker products; we provide the explicit formulas in Appendix~\ref{sec:appendix-local-inverse}.


\section{Numerical results}
\label{sec:results}

In this section, we present a series of numerical experiments to evaluate the performance of the proposed local
p-multigrid solvers and their effectiveness as smoothers in a global geometric multigrid hierarchy for the Stokes
equations. We begin by validating the local solvers on a single patch, allowing us to isolate their convergence
properties from the global multigrid components and assess their sensitivity to mesh distortion and coefficient
variation. Subsequently, we incorporate these local solvers into a full geometric multigrid framework and investigate
the overall solver efficiency for various problem configurations.

All tests are conducted using a homogeneous right-hand side ($f=0$) and a random initial guess, which provides a
representative test for the smoothing properties across all frequency modes. We note that the current experimental
implementation is based on sparse matrix data structures to facilitate development (and consequently, debugging). While
this choice allows for rapid prototyping of different local solver variants, it does not reflect the peak performance
achievable with a fully optimized matrix-free implementation. Consequently, our performance metrics focus primarily on
iteration counts and relative convergence rates rather than absolute runtimes, which would be significantly lower in a
production-ready matrix-free environment.

\subsection{Numerical validation on a single patch}
\label{sec:local_solve_validation}

Before integrating the local p-multigrid solver into the global smoothing procedure, we first evaluate its performance
on a single vertex patch. This validation section follows the methodology presented in~\cite{wichrowski2025local} for
the scalar Laplace problem, which we adapt here for the Stokes equations. For these tests, we run the FGMRES solver
using a p-multigrid V-cycle as a preconditioner until the relative residual is reduced below a tolerance of $10^{-8}$.
We emphasize that such a stringent tolerance is not necessary when the p-multigrid is used as a patch smoother in the
global algorithm: a much looser tolerance typically suffices and is computationally cheaper. We nevertheless enforce
$10^{-8}$ here to increase the number of FGMRES iterations, which makes convergence behavior easier to inspect and
provides more detailed diagnostics for validation. These tests are not intended to validate the local solver as a
standalone production solver; rather, they isolate its behaviour from other multigrid components and assess only the
contributions proposed in this paper.

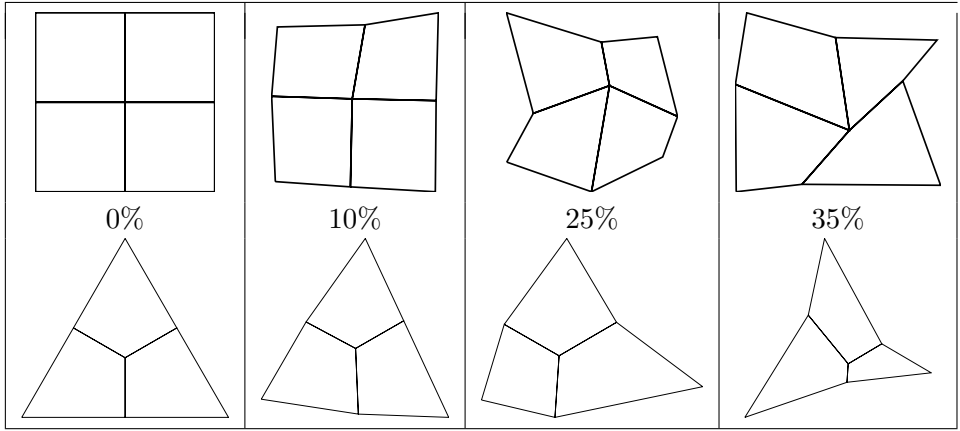
\begin{figure}[htbp]
  \centering
\begin{tabular}{|c|c|c|c|}
  \hline
      &      &             \\[-0.9em]
  \begin{tikzpicture}[scale=0.8, every node/.style={scale=0.8}]
    \begin{axis}[
        axis equal image,
        hide axis,
        width=0.32\columnwidth,
        enlargelimits=false
      ]
      \addplot[black,thick,mark=none] table {figures/patch_def_data/patch00_2D_2.gp};
    \end{axis}
  \end{tikzpicture}
      &
  \begin{tikzpicture}[scale=0.8, every node/.style={scale=0.8}]
    \begin{axis}[
        axis equal image,
        hide axis,
        width=0.32\columnwidth,
        enlargelimits=false
      ]
      \addplot[black,thick,mark=none] table {figures/patch_def_data/patch01_2D.gp};
    \end{axis}
  \end{tikzpicture}
      &
  \begin{tikzpicture}[scale=0.8, every node/.style={scale=0.8}]
    \begin{axis}[
        axis equal image,
        hide axis,
        width=0.32\columnwidth,
        enlargelimits=false
      ]
      \addplot[black,thick,mark=none] table {figures/patch_def_data/patch025_2D_1.gp};
    \end{axis}
  \end{tikzpicture}
      &
  \begin{tikzpicture}[scale=0.8, every node/.style={scale=0.8}]
    \begin{axis}[
        axis equal image,
        hide axis,
        width=0.32\columnwidth,
        enlargelimits=false
      ]
      \addplot[black,thick,mark=none] table {figures/patch_def_data/patch035_2D.gp};
    \end{axis}
  \end{tikzpicture}
  \\
  0\% & 10\% & 25\% & 35\%
  \\
  \begin{tikzpicture}[scale=0.8, every node/.style={scale=0.8}]
    \begin{axis}[
        axis equal image,
        hide axis,
        width=0.32\columnwidth,
        enlargelimits=false
      ]
      \addplot[black,mark=none] table {figures/patch_def_data/simplex00_2D_2.gp};
    \end{axis}
  \end{tikzpicture}
      &
  \begin{tikzpicture}[scale=0.8, every node/.style={scale=0.8}]
    \begin{axis}[
        axis equal image,
        hide axis,
        width=0.32\columnwidth,
        enlargelimits=false
      ]
      \addplot[black,mark=none] table {figures/patch_def_data/simplex01_2D.gp};
    \end{axis}
  \end{tikzpicture}
      &
  \begin{tikzpicture}[scale=0.8, every node/.style={scale=0.8}]
    \begin{axis}[
        axis equal image,
        hide axis,
        width=0.32\columnwidth,
        enlargelimits=false
      ]
      \addplot[black,mark=none] table {figures/patch_def_data/simplex025_2D.gp};
    \end{axis}
  \end{tikzpicture}
      &
  \begin{tikzpicture}[scale=0.8, every node/.style={scale=0.8}]
    \begin{axis}[
        axis equal image,
        hide axis,
        width=0.32\columnwidth,
        enlargelimits=false
      ]
      \addplot[black,mark=none] table {figures/patch_def_data/simplex035_2D.gp};
    \end{axis}
  \end{tikzpicture}
  \\
  \hline
\end{tabular}
  \caption{Illustration of the vertex-based patches used in our 2D experiments. The top row shows structured (Cartesian) patches, while the bottom row depicts unstructured (simplicial) patches. The distortion level $\delta$ represents the maximum displacement of each interior vertex relative to the local mesh size $h$.
    \label{fig:deformed_patches}
  }
\end{figure}

We evaluate the performance of various local solver variants on a single vertex patch consisting of $2^d$ cells in
spatial dimensions $d=2$ and $d=3$. Following~\cite{wichrowski2025local}, we also test simplicial (unstructured)
patches comprising 3 quadrilaterals in 2D or 4 hexahedra in 3D; representative geometries are shown in
Figure~\ref{fig:deformed_patches}. Our investigation starts with the Stokes problem on an undeformed Cartesian patch as
a baseline and progressively increases difficulty to probe the robustness and limitations of the proposed methods.

The first level of complexity involves geometric distortion introduced by randomly displacing each vertex. As
in~\cite{wichrowski2025local}, each interior vertex is shifted by a vector of length $\delta h$, with its direction
sampled uniformly on the unit sphere $S^{d-1}$. For larger distortions, generated patches may contain degenerate
elements with concave or inverted corners. For 2D Cartesian patches, degeneracy occurs when the distortion exceeds
$\delta = \frac{1}{2}\frac{\sqrt{2}}{2} \approx 0.3535$, while in 3D the threshold is $\delta =
  \frac{1}{2}\frac{\sqrt{3}}{3} \approx 0.2887$. Accordingly, we test distortions up to 35\% in 2D and 30\% in 3D. We
note that simplicial patches exhibit higher sensitivity to distortion because the initial angles already exceed
$90^\circ$, causing some cells to approach triangular shapes with angles near $180^\circ$.

Secondly, we investigate the influence of heterogeneous material properties, focusing on discontinuous viscosity jumps
where one cell carries a large contrast relative to others. Finally, we combine geometric deformation and strong
coefficient variation to evaluate the solvers under extreme conditions, considering both stationary and non-stationary
flow regimes. For each scenario, we report FGMRES iteration counts (to a tolerance of $10^{-8}$) across a range of
polynomial degrees $p$. Following the findings in~\cite{wichrowski2025local}, where 2D and 3D results were shown to be
qualitatively similar, we focus most experiments on the 2D case to save computational resources and present 3D results
only for selected configurations.

We perform 20 independent random realizations when testing geometric distortion and report the average number of
iterations over all runs. Each local solve is allowed a maximum of 150 iterations; if a run does not reach the
convergence tolerance within this limit, it is declared non-convergent. Such runs would be assigned a large sentinel
iteration count to make non-convergence immediately apparent.

\paragraph{Convergence on Cartesian Patches}
We begin by evaluating the baseline performance of the various local solvers on undeformed Cartesian patches with a
constant coefficient. In this setting, the velocity block $A_j$ can be inverted exactly using the Fast Diagonalization,
providing a reference case for evaluating the impact of different approximations for the velocity and Schur complement
blocks.

We compare five different variants: (i) the Braess--Sarazin smoother using a diagonal approximation $\tilde{A}_j^{-1}$
and an inner Richardson iteration for the Schur complement; (ii) Braess--Sarazin using the same $\tilde{A}_j^{-1}$ but
with a single CG iteration for the Schur system; (iii) a block GMRES solver using the exact velocity inverse and a
pressure mass matrix $\mathcal{M}_p^{-1}$ for the Schur contribution; (iv) block GMRES where the velocity block is
instead approximated by a single $p$-multigrid V-cycle; and (v) the blockwise elimination approach using the exact
velocity inverse and an iterative CG solve for the Schur complement. Table~\ref{tab:single_patch_cartesian} summarizes
the number of FGMRES iterations required for convergence as a function of the polynomial degree $p$.

We observe that the Braess--Sarazin variants demonstrate remarkable robustness: even with the simple diagonal
approximation for the velocity block, they consistently outperform the block-based solvers. For instance, with $m=2$
smoothing steps, the Braess--Sarazin iteration counts are nearly factor of two lower than those of the block GMRES and
blockwise elimination methods, even when the latter are equipped with the exact inverse $A^{-1}$. Under these ideal
Cartesian conditions, the CG-based and Richardson-based inner Schur solvers for Braess--Sarazin yield virtually
identical iteration counts.

In contrast, the performance of the block GMRES solver is highly sensitive to the quality of the velocity inverse.
While the iteration counts are stable when $A^{-1}$ is inverted exactly, substituting it with a single $p$-multigrid
V-cycle (variant iv) leads to a significant degradation in performance, with the iteration count rising sharply as the
polynomial degree $p$ increases.

We note that the iteration counts exhibit characteristic non-monotonic behavior, with visible drops at $p=4$ and $p=8$.
This is an artifact of the $p$-multigrid level structure used in the hierarchy, an effect previously analyzed for the
scalar case in~\cite{wichrowski2025local}. We build $p$-levels using the nested sequence $1, 3, 7, 15$, etc. (i.e.,
$p_{k+1} = 2p_k + 1$). Degrees that are just above these nodes benefit from a more complete multigrid hierarchy; for
example, for $p=7$ the hierarchy is simply $\{1,3, 7\}$, whereas for $p=8$ it becomes $\{1, 3, 7, 8\}$, effectively
yielding a method closer to performing two smoothing steps on the finest level.

\begin{table}[htbp]
  \centering
  \caption{Number of iterations for solving the local Stokes problem on a 2D Cartesian patch. Comparison of Braess-Sarazin with two variants of Schur solver: Richardson and CG, Block GMRES with exact velocity inverse and p-multigrid velocity inverse, and blockwise elimination with exact velocity inverse. }
  \label{tab:single_patch_cartesian}
  \begin{tabular}{lccc|ccccccc}
    \toprule
    Method         & $\tilde{A}^{-1}$                        & $\tilde{S}^{-1}$           & $m$ & \multicolumn{7}{c}{Polynomial Degree $p$}                               \\
                   &                                         &                            &     & 2                                         & 3  & 4  & 5  & 7  & 8  & 11 \\
    \midrule
    \multirow{4}{*}{Braess-Sarazin}
                   & \multirow{4}{*}{$\omega \diag(A)^{-1}$}
                   & \multirow{2}{*}{$\diag(\hat{S}^{-1})$}  & 1                          & 6   & 9                                         & 8  & 9  & 13 & 10 & 13      \\
                   &                                         &                            & 2   & 4                                         & 6  & 6  & 6  & 9  & 7  & 9  \\
    \cline{3-11}
                   &                                         & \multirow{2}{*}{1 PCG it.} & 1   & 6                                         & 10 & 8  & 9  & 13 & 10 & 13 \\
                   &                                         &                            & 2   & 4                                         & 6  & 6  & 6  & 9  & 7  & 9  \\
    \cline{1-11}
    \multirow{3}{*}{Block GMRES}
                   & \multirow{2}{*}{p-MG V-cycle}
                   & \multirow{2}{*}{$M_p^{-1}$}             & 1                          & 15  & 28                                        & 33 & 38 & 52 & 45 & 57      \\
                   &                                         &                            & 2   & 13                                        & 24 & 27 & 31 & 39 & 37 & 44 \\
    \cline{2-11}
    \myrowcolor
                   & $ A^{-1}$                               & $M_p^{-1}$                 & --  & 9                                         & 14 & 16 & 17 & 18 & 19 & 20 \\
    \cline{1-11}
    Blockwise Elm. & $A^{-1}$                                & CG                         & --  & 8                                         & 14 & 16 & 18 & 21 & 22 & 22 \\
    \bottomrule
  \end{tabular}
\end{table}

\begin{figure}[htbp]
  \begin{tikzpicture}
    \begin{scope}
      \begin{axis}[
          paperplot,
          xlabel={distortion},
          ylabel={Iteration count},
          xmin=0.0, xmax=0.35,
          xtick={0.05,0.1,0.15,0.2,0.25,0.3,0.35},
          ymode=log,
          log basis y=8,
          ytick={2,4,8,16,32,64,128},
          yticklabels={2,4,8,16,32,64,128},
          minor y tick num=3,
          xmajorgrids=true,
          ymajorgrids=true,
          yminorgrids=true,
          ymin=4,
          ymax=80,
          legend pos=south east,
          legend columns=3,
          axis lines=box,
          clip=true,
        ]
        \node[anchor=north east,inner sep=5pt,fill=white,fill opacity=0.9,text opacity=1,draw=black,rounded
          corners=2pt] at (rel axis cs:0.2,0.92) {2D};

        \addplot+ [p2, solid] table [col sep=comma, x=distortion, y=avg] {results/tmp/struct_2d_block_p2.csv}; \addlegendentry{$p=2$}
        \addplot+ [p3, solid] table [col sep=comma, x=distortion, y=avg] {results/tmp/struct_2d_block_p3.csv}; \addlegendentry{$p=3$}
        \addplot+ [p7, solid] table [col sep=comma, x=distortion, y=avg] {results/tmp/struct_2d_block_p7.csv}; \addlegendentry{$p=7$}

        \addplot+ [p2, dashed, forget plot] table [col sep=comma, x=distortion, y=avg] {results/tmp/struct_2d_exact_p2.csv};
        \addplot+ [p3, dashed, forget plot] table [col sep=comma, x=distortion, y=avg] {results/tmp/struct_2d_exact_p3.csv};
        \addplot+ [p7, dashed, forget plot] table [col sep=comma, x=distortion, y=avg] {results/tmp/struct_2d_exact_p7.csv};

      \end{axis}
    \end{scope}

    \begin{scope}[xshift=0.52\columnwidth]
      \begin{axis}[
          paperplot,
          xlabel={distortion},
          ylabel={Iteration count},
          xmin=0.0, xmax=0.3,
          xtick={0.05,0.1,0.15,0.2,0.25,0.3,0.35},
          ymode=log,
          log basis y=8,
          ytick={2,4,8,16,32,64,128},
          yticklabels={2,4,8,16,32,64,128},
          minor y tick num=3,
          xmajorgrids=true,
          ymajorgrids=true,
          yminorgrids=true,
          ymin=4,
          legend pos=south east,
          legend columns=3,
        ]
        \node[anchor=north east,inner sep=5pt,fill=white,fill opacity=0.9,text opacity=1,draw=black,rounded
          corners=2pt] at (rel axis cs:0.2,0.92) {3D};

        \addplot+ [p2, solid] table [col sep=comma, x=distortion, y=avg] {results/tmp/struct_3d_block_p2.csv}; \addlegendentry{$p=2$}
        \addplot+ [p3, solid] table [col sep=comma, x=distortion, y=avg] {results/tmp/struct_3d_block_p3.csv}; \addlegendentry{$p=3$}
        \addplot+ [p7, solid] table [col sep=comma, x=distortion, y=avg] {results/tmp/struct_3d_block_p7.csv}; \addlegendentry{$p=7$}

        \addplot+ [p2, dashed, forget plot] table [col sep=comma, x=distortion, y=avg] {results/tmp/struct_3d_exact_p2.csv};
        \addplot+ [p3, dashed, forget plot] table [col sep=comma, x=distortion, y=avg] {results/tmp/struct_3d_exact_p3.csv};
        \addplot+ [p7, dashed, forget plot] table [col sep=comma, x=distortion, y=avg] {results/tmp/struct_3d_exact_p7.csv};

      \end{axis}
    \end{scope}
  \end{tikzpicture}
  \caption{Number of FGMRES iterations (logarithmic y-axis) vs.\ mesh distortion for local block solver on structured patches in 2D (left) and 3D (right). Solid lines: velocity block $A$ is approximated via a single $p$-multigrid V-cycle; dashed lines: exact velocity inverse; dotted lines: reference.}
  \label{fig:baseline_iter_vs_distortion_block}
\end{figure}
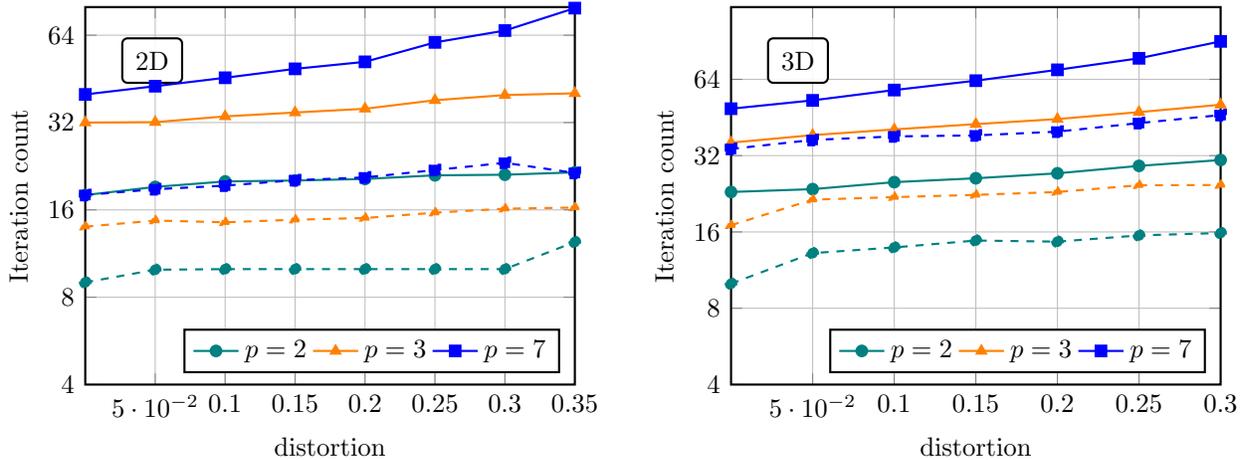

\paragraph{Resilience to Geometric Distortion}
Next, we examine the solver's performance under geometric perturbations. We apply a random, nonlinear distortion to the
patch geometry controlled by a \emph{distortion} parameter $\delta$, with zero corresponding to an undeformed patch.
Figure~\ref{fig:baseline_iter_vs_distortion_block} presents the results for the block GMRES solvers used as a baseline
for comparison. Here, the local Stokes problem is solved using a block preconditioner where the velocity block $A$ is
approximated by either a single $p$-multigrid V-cycle (solid lines) or its exact inverse (dashed lines). The left and
right panels show results for 2D and 3D, respectively. For small distortions, the iteration counts are stable across
all polynomial degrees. However, as the distortion increases, the performance depends on the spatial dimension and the
polynomial degree. In 2D, the iteration count grows moderately, and the gap between the $p$-multigrid approximation and
the exact inverse remains small. In 3D, the block solver is significantly more sensitive: for higher polynomial degrees
such as $p=7$, the iteration counts rise faster with distortion.

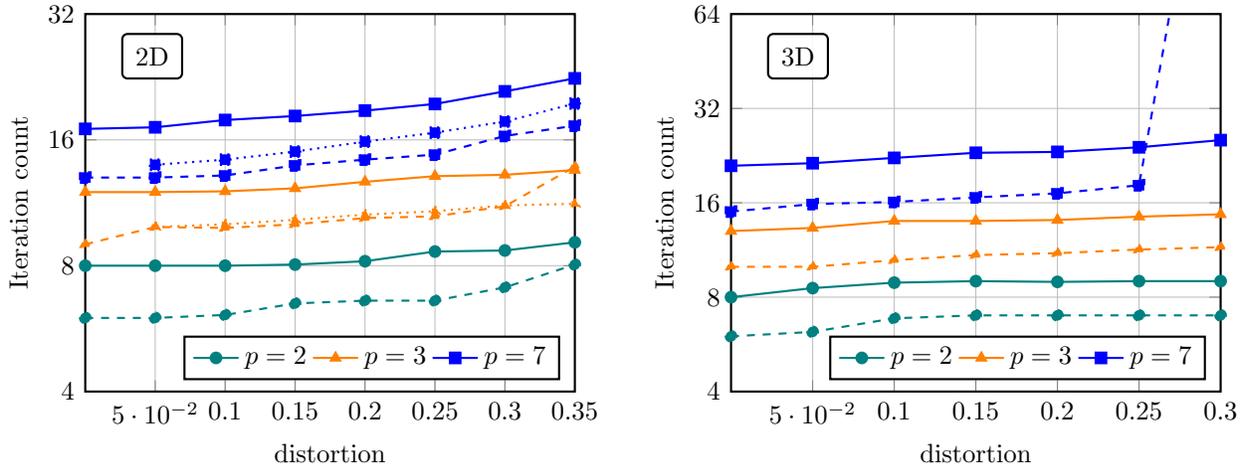
\begin{figure}[htbp]
  \begin{tikzpicture}
    \begin{scope}
      \begin{axis}[
          paperplot,
          xlabel={distortion},
          ylabel={Iteration count},
          xmin=0.0, xmax=0.35,
          xtick={0.05,0.1,0.15,0.2,0.25,0.3,0.35},
          ymode=log,
          log basis y=8,
          ytick={2,4,8,16,32,64,128},
          yticklabels={2,4,8,16,32,64,128},
          minor y tick num=3,
          xmajorgrids=true,
          ymajorgrids=true,
          yminorgrids=true,
          ymin=4,
          ymax=32,
          legend pos=south east,
          axis lines=box,
          clip=true,
          legend columns=3,
        ]
        \node[anchor=north east,inner sep=5pt,fill=white,fill opacity=0.9,text opacity=1,draw=black,rounded
          corners=2pt] at (rel axis cs:0.2,0.95) {2D};

        \addplot+ [p2, solid] table [col sep=comma, x=distortion, y=avg] {results/tmp/struct_2d_braess_p2.csv}; \addlegendentry{$p=2$}
        \addplot+ [p3, solid] table [col sep=comma, x=distortion, y=avg] {results/tmp/struct_2d_braess_p3.csv}; \addlegendentry{$p=3$}
        \addplot+ [p7, solid] table [col sep=comma, x=distortion, y=avg] {results/tmp/struct_2d_braess_p7.csv}; \addlegendentry{$p=7$}

        \addplot+ [p2, dashed, forget plot] table [col sep=comma, x=distortion, y=avg] {results/tmp/simplex_2d_braess_p2.csv};
        \addplot+ [p3, dashed, forget plot] table [col sep=comma, x=distortion, y=avg] {results/tmp/simplex_2d_braess_p3.csv};
        \addplot+ [p7, dashed, forget plot] table [col sep=comma, x=distortion, y=avg] {results/tmp/simplex_2d_braess_p7.csv};

        \addplot+ [p3, dotted, thick, forget plot] table [col sep=comma, x=distortion, y=avg] {results/tmp/poisson_2d_p3.csv};
        \addplot+ [p7, dotted, thick, forget plot] table [col sep=comma, x=distortion, y=avg] {results/tmp/poisson_2d_p7.csv};

      \end{axis}
    \end{scope}

    \begin{scope}[xshift=0.52\columnwidth]
      \begin{axis}[
          paperplot,
          xlabel={distortion},
          ylabel={Iteration count},
          xmin=0.0, xmax=0.3,
          xtick={0.05,0.1,0.15,0.2,0.25,0.3,0.35},
          ymode=log,
          log basis y=8,
          ytick={2,4,8,16,32,64,128},
          yticklabels={2,4,8,16,32,64,128},
          minor y tick num=3,
          xmajorgrids=true,
          ymajorgrids=true,
          yminorgrids=true,
          ymin=4,
          ymax=64,
          legend pos=south east,
          legend columns=3,
        ]
        \node[anchor=north east,inner sep=5pt,fill=white,fill opacity=0.9,text opacity=1,draw=black,rounded
          corners=2pt] at (rel axis cs:0.2,0.95) {3D};

        \addplot+ [p2, solid] table [col sep=comma, x=distortion, y=avg] {results/tmp/struct_3d_braess_p2.csv}; \addlegendentry{$p=2$}
        \addplot+ [p3, solid] table [col sep=comma, x=distortion, y=avg] {results/tmp/struct_3d_braess_p3.csv}; \addlegendentry{$p=3$}
        \addplot+ [p7, solid] table [col sep=comma, x=distortion, y=avg] {results/tmp/struct_3d_braess_p7.csv}; \addlegendentry{$p=7$}

        \addplot+ [p2, dashed, forget plot] table [col sep=comma, x=distortion, y=avg] {results/tmp/simplex_3d_braess_p2.csv};
        \addplot+ [p3, dashed, forget plot] table [col sep=comma, x=distortion, y=avg] {results/tmp/simplex_3d_braess_p3.csv};
        \addplot+ [p7, dashed, forget plot] table [col sep=comma, x=distortion, y=avg] {results/tmp/simplex_3d_braess_p7.csv};

      \end{axis}
    \end{scope}
  \end{tikzpicture}
  \caption{Number of FGMRES iterations (logarithmic y-axis) vs.\ mesh distortion for Braess--Sarazin smoother on structured patches. Both variants use one smoothing step per level. Solid lines: structured patch, dashed lines: simplicial patch, dotted lines: reference Poisson problem iteration counts (available for $p=3$ and $p=7$). Left subplot: 2D; right subplot: 3D. For $p=7$ the solver did not converge in some runs with $\delta=0.3$ on simplicial patches. Dotted lines: reference Poisson problem iteration counts (available for $p=3$ and $p=7$ in 2D)}
  \label{fig:iter_vs_distortion_structured}
\end{figure}

We now contrast these baseline results with the Braess--Sarazin smoother, which is our primary focus. In
Figure~\ref{fig:iter_vs_distortion_structured}, we show the iteration counts (again on a logarithmic scale) for the
p-multigrid solver using Braess--Sarazin smoothing. Figure~\ref{fig:iter_vs_distortion_structured} shows the iteration
counts for the Braess--Sarazin smoother using one smoothing step per level ($m=1$). We compare structured patches
(solid lines) with unstructured simplex patches (dashed lines). In 2D (left panel), the Braess--Sarazin smoother
demonstrates remarkable robustness, with iteration counts remaining below 20 even for 35\% distortion, showing almost
no dependence on the polynomial degree. The only exception is the highest polynomial degree $p=7$ on simplicial
patches, where the solver failed to converge for some realizations with large distortion $\delta=0.3$. In 3D (right
panel), while iteration counts are slightly higher, they remain stable and significantly more efficient than the block
solvers from Figure~\ref{fig:baseline_iter_vs_distortion_block}. Notably, the Braess--Sarazin smoother handles both
structured and simplicial patches with similar efficiency, and the performance degradation observed for high $p$ in 3D
is much less pronounced here. This suggests that the saddle-point coupling in the Braess--Sarazin approach is better
suited for distorted patches than the decoupled block-preconditioning approach. We also compare the results with those
for the scalar Poisson problem on structured patches (dotted lines) from~\cite{wichrowski2025local}, which shows
similar trends while requiring slightly fewer iterations.

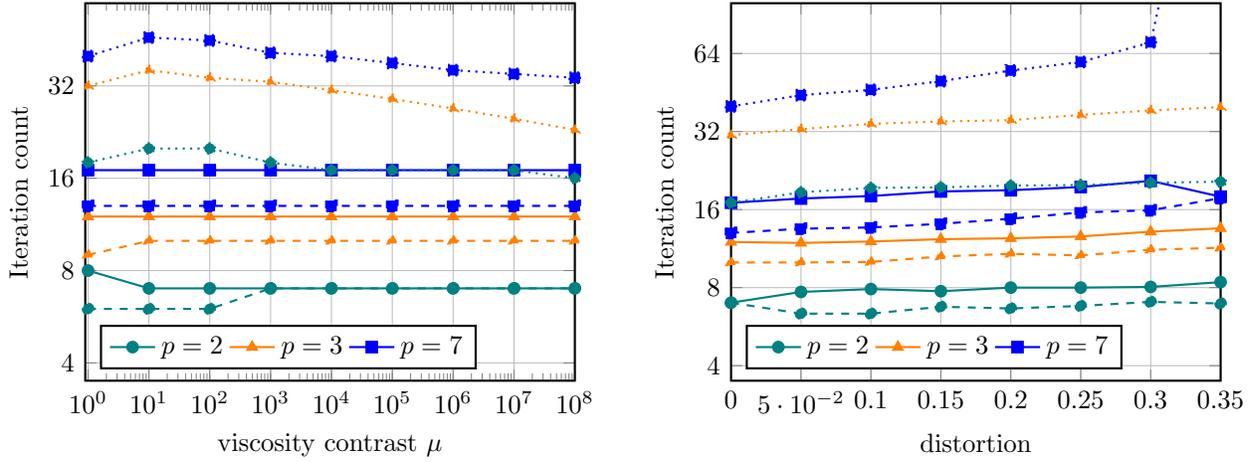
\begin{figure}[htbp]
  \centering
  \begin{tikzpicture}
    \begin{scope}
      \begin{axis}[
          paperplot,
          xlabel={viscosity contrast $\mu$},
          ylabel={Iteration count},
          xmin=.9,
          xmax=1e8,
          xmode=log,
          log basis x=10,
          ymode=log,
          log basis y=8,
          ytick={2,4,8,16,32,64,128},
          yticklabels={2,4,8,16,32,64,128},
          minor y tick num=4,
          xmajorgrids=true,
          ymajorgrids=true,
          yminorgrids=true,
          ymin=3.5,
          legend pos=south west,
          legend columns=3,
        ]
        \addplot+ [p2, solid] table [col sep=comma, x=mu, y=avg] {results/tmp/struct_2d_braess_jump_p2.csv}; \addlegendentry{$p=2$}
        \addplot+ [p3, solid] table [col sep=comma, x=mu, y=avg] {results/tmp/struct_2d_braess_jump_p3.csv}; \addlegendentry{$p=3$}
        \addplot+ [p7, solid] table [col sep=comma, x=mu, y=avg] {results/tmp/struct_2d_braess_jump_p7.csv}; \addlegendentry{$p=7$}

        \addplot+ [p2, dashed, forget plot] table [col sep=comma, x=mu, y=avg] {results/tmp/simplex_2d_braess_jump_p2.csv};
        \addplot+ [p3, dashed, forget plot] table [col sep=comma, x=mu, y=avg] {results/tmp/simplex_2d_braess_jump_p3.csv};
        \addplot+ [p7, dashed, forget plot] table [col sep=comma, x=mu, y=avg] {results/tmp/simplex_2d_braess_jump_p7.csv};

        \addplot+ [p2, dotted, thick, forget plot] table [col sep=comma, x=mu, y=avg] {results/tmp/struct_2d_block_jump_p2.csv};
        \addplot+ [p3, dotted, thick, forget plot] table [col sep=comma, x=mu, y=avg] {results/tmp/struct_2d_block_jump_p3.csv};
        \addplot+ [p7, dotted, thick, forget plot] table [col sep=comma, x=mu, y=avg] {results/tmp/struct_2d_block_jump_p7.csv};

      \end{axis}
    \end{scope}

    \begin{scope}[xshift=0.52\columnwidth]
      \begin{axis}[
          paperplot,
          xlabel={distortion},
          ylabel={Iteration count},
          xmin=0,
          xmax=0.35,
          ymode=log,
          log basis y=8,
          ytick={2,4,8,16,32,64,128},
          yticklabels={2,4,8,16,32,64,128},
          minor y tick num=3,
          xmajorgrids=true,
          ymajorgrids=true,
          yminorgrids=true,
          ymin=3.5,
          ymax=100,
          legend pos=south west,
          legend columns=3,
        ]
        \addplot+ [p2, solid] table [col sep=comma, x=distortion, y=avg] {results/tmp/struct_2d_braess_dist_mu4_p2.csv}; \addlegendentry{$p=2$}
        \addplot+ [p3, solid] table [col sep=comma, x=distortion, y=avg] {results/tmp/struct_2d_braess_dist_mu4_p3.csv}; \addlegendentry{$p=3$}
        \addplot+ [p7, solid] table [col sep=comma, x=distortion, y=avg] {results/tmp/struct_2d_braess_dist_mu4_p7.csv}; \addlegendentry{$p=7$}

        \addplot+ [p2, dashed, forget plot] table [col sep=comma, x=distortion, y=avg] {results/tmp/simplex_2d_braess_dist_mu4_p2.csv};
        \addplot+ [p3, dashed, forget plot] table [col sep=comma, x=distortion, y=avg] {results/tmp/simplex_2d_braess_dist_mu4_p3.csv};
        \addplot+ [p7, dashed, forget plot] table [col sep=comma, x=distortion, y=avg] {results/tmp/simplex_2d_braess_dist_mu4_p7.csv};

        \addplot+ [p2, dotted, thick, forget plot] table [col sep=comma, x=distortion, y=avg] {results/tmp/struct_2d_block_dist_mu4_p2.csv};
        \addplot+ [p3, dotted, thick, forget plot] table [col sep=comma, x=distortion, y=avg] {results/tmp/struct_2d_block_dist_mu4_p3.csv};
        \addplot+ [p7, dotted, thick, forget plot] table [col sep=comma, x=distortion, y=avg] {results/tmp/struct_2d_block_dist_mu4_p7.csv};

      \end{axis}
    \end{scope}
  \end{tikzpicture}
  \caption{Solver iterations (logarithmic y-axis) vs.\ coefficient jump and mesh distortion (Stationary case, $\rho=0$).
    Left panel: FGMRES iterations on an undeformed 2D patch as the coefficient jump $\mu$ varies.
    Right panel: FGMRES iterations on a 2D patch with a fixed coefficient jump $\mu=10^4$ as the mesh distortion increases.
    Solid lines: structured patches; dashed lines: Simplex patches; dotted lines: block solver.
    In both panels, we use p-multigrid with Braess--Sarazin smoother and one smoothing step per level. For $p=7$ the block solver did not converge in some runs with $\delta=0.35$.}
  \label{fig:iter_vs_jump_stationary}
\end{figure}

\paragraph{Viscosity Jumps and Non-Zero Density }
We now investigate robustness to variations in the material coefficient, focusing on a challenging discontinuous jump.
In this test, one cell in the patch is assigned a coefficient $\mu$ while all other cells have coefficient 1.
Figure~\ref{fig:iter_vs_jump_stationary} presents two related experiments in 2D using the Braess--Sarazin smoother. The
left panel plots the number of iterations (shown on a logarithmic y-axis) against the magnitude of the coefficient jump
$\mu$ on an undeformed patch. The Braess--Sarazin smoother (solid and dashed lines for Cartesian and simplex patches,
respectively) demonstrates excellent robustness, with iteration counts showing only a minimal increase as the jump
magnitude grows from $10^0$ to $10^8$. In this configuration, the block GMRES solver (dotted lines), which utilizes the
pressure mass matrix inverse $M_p^{-1}$ as a preconditoner for the Schur complement, also maintains stable performance.
While it requires roughly twice as many iterations as the Braess--Sarazin approach, it does not exhibit any severe
degradation as the jump magnitude increases. These results suggest that Braess--Sarazin remains the more efficient and
consistent choice across varying degrees of numerical difficulty, providing significant acceleration compared to the
simpler decoupled block-preconditioning approach.

The right panel of Figure~\ref{fig:iter_vs_jump_stationary} examines the combined challenge of geometric distortion and
a large coefficient jump. It shows the iteration counts for a 2D patch with a fixed jump $\mu=10^4$ as the distortion
increases. By comparing these results with the constant-coefficient case
(Figure~\ref{fig:iter_vs_distortion_structured}), we observe that the presence of the jump has a negligible effect on
the Braess--Sarazin solver's sensitivity to distortion. However, the block solver (dotted lines) did not converge for
$p=7$ when distortion reached 35\%. This level of robustness is critical for practical applications where complex
geometries and heterogeneous materials are often encountered together.

\begin{figure}[htbp]
  \begin{tikzpicture}
    \begin{scope}
      \begin{axis}[
          paperplot,
          xlabel={viscosity contrast $\mu$},
          ylabel={Iteration count},
          xmin=.9,
          xmax=1e8,
          xmode=log,
          ymode=log,
          log basis x=10,
          log basis y=8,
          ytick={2,4,8,16,32,64,128},
          yticklabels={2,4,8,16,32,64,128},
          minor y tick num=3,
          xmajorgrids=true,
          ymajorgrids=true,
          yminorgrids=true,
          ymin=3.5,
          legend pos=south west,
          legend columns=3,
        ]
        \addplot+ [p2, solid] table [col sep=comma, x=mu, y=avg] {results/tmp/struct_2d_braess_jump_nonstat_p2.csv}; \addlegendentry{$p=2$}
        \addplot+ [p3, solid] table [col sep=comma, x=mu, y=avg] {results/tmp/struct_2d_braess_jump_nonstat_p3.csv}; \addlegendentry{$p=3$}
        \addplot+ [p7, solid] table [col sep=comma, x=mu, y=avg] {results/tmp/struct_2d_braess_jump_nonstat_p7.csv}; \addlegendentry{$p=7$}

        \addplot+ [p2, dotted, thick, forget plot] table [col sep=comma, x=mu, y=avg] {results/tmp/struct_2d_block_jump_nonstat_p2.csv};
        \addplot+ [p3, dotted, thick, forget plot] table [col sep=comma, x=mu, y=avg] {results/tmp/struct_2d_block_jump_nonstat_p3.csv};
        \addplot+ [p7, dotted, thick, forget plot] table [col sep=comma, x=mu, y=avg] {results/tmp/struct_2d_block_jump_nonstat_p7.csv};
      \end{axis}
    \end{scope}

    \begin{scope}[xshift=0.52\columnwidth]
      \begin{axis}[
          paperplot,
          xlabel={distortion},
          ylabel={Iteration count},
          xmin=0,
          xmax=0.35,
          ymode=log,
          log basis y=8,
          ytick={2,4,8,16,32,64,128},
          yticklabels={2,4,8,16,32,64,128},
          minor y tick num=3,
          xmajorgrids=true,
          ymajorgrids=true,
          yminorgrids=true,
          ymin=3.5,
          ymax=128,
          legend pos=south west,
          legend columns=3,
        ]

        \addplot+ [p2, solid] table [col sep=comma, x=distortion, y=avg]
          {results/tmp/struct_2d_braess_dist_nonstat_mu4_p2.csv};
        \addlegendentry{$p=2$}

        \addplot+ [p3, solid] table [col sep=comma, x=distortion, y=avg]
          {results/tmp/struct_2d_braess_dist_nonstat_mu4_p3.csv};
        \addlegendentry{$p=3$}

        \addplot+ [p7, solid] table [col sep=comma, x=distortion, y=avg]
          {results/tmp/struct_2d_braess_dist_nonstat_mu4_p7.csv};
        \addlegendentry{$p=7$}

        \addplot+ [p2, dotted, thick, forget plot] table [col sep=comma, x=distortion, y=avg]
          {results/tmp/struct_2d_block_dist_nonstat_mu4_p2.csv};
        \addplot+ [p3, dotted, thick, forget plot] table [col sep=comma, x=distortion, y=avg]
          {results/tmp/struct_2d_block_dist_nonstat_mu4_p3.csv};
        \addplot+ [p7, dotted, thick, forget plot] table [col sep=comma, x=distortion, y=avg]
          {results/tmp/struct_2d_block_dist_nonstat_mu4_p7.csv};

      \end{axis}
    \end{scope}
  \end{tikzpicture}
  \caption{Solver iterations (logarithmic y-axis) vs.\ coefficient jump and mesh distortion for non-stationary Stokes ($\rho=1$).
    Left panel: iterations on an undeformed 2D structured patch as viscosity jump $\mu$ varies.
    Right panel: iterations on a 2D structured patch with fixed jump $\mu=10^4$ and increasing distortion.
    Solid lines: Braess--Sarazin smoother; dotted lines: block solver. Both variants use p-multigrid for the velocity block. For $p=7$ the block solver did not converge in some runs with $\delta\geq 0.3$.}
  \label{fig:iter_vs_jump_nonstationary}
\end{figure}
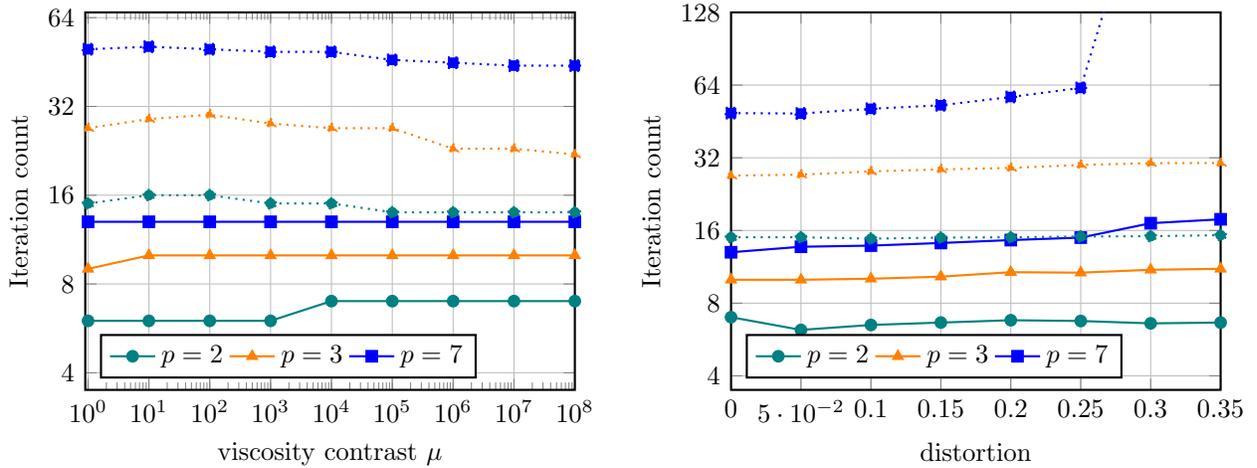

Finally, we evaluate the solver's robustness for the non-stationary Stokes equations by setting the density $\rho=1$
(Figure~\ref{fig:iter_vs_jump_nonstationary}). The results for the Braess--Sarazin smoother (solid lines) are
qualitatively similar to the stationary case, with excellent robustness to both viscosity jumps and mesh distortion. In
this non-stationary regime, the block GMRES solver (dotted lines) also demonstrates stable performance across the
entire range of viscosity contrasts, although it fails to converge for high polynomial degrees ($p=7$) when the patch
distortion exceeds 30\%. While it requires approximately twice as many iterations as the Braess--Sarazin approach, it
avoids the typical convergence issues associated with strong coefficient jumps on less distorted meshes. Notably, the
iteration counts for the block solver are nearly identical to those observed in the stationary case, indicating that
the addition of the mass term in the velocity block does not significantly alter the solver's asymptotic behavior with
respect to the jump magnitude in this configuration. As in the stationary case, the block solver baseline uses
$M_p^{-1}$ to approximate the Schur complement inverse; although there is room for further improvement through more
advanced Schur preconditioners, Braess--Sarazin remains the more efficient choice for rapid and consistent convergence.

\paragraph{Summary}
Our numerical study on single patches demonstrates that the local p-multigrid solver using the Braess--Sarazin smoother
is a highly robust and versatile choice. It maintains stable performance across a wide range of mesh distortions,
spatial dimensions, and extreme coefficient jumps, making it well-suited for the complex conditions encountered in
practical Stokes flow simulations. While its iteration counts show some sensitivity to the polynomial degree $p$ and
the mesh distortion, the degradation is much more graceful than the decoupled block-preconditioning approaches. The
performance gap is clearly visible in the logarithmic plots shown in Figures~\ref{fig:iter_vs_jump_stationary} and
\ref{fig:iter_vs_jump_nonstationary}.

In contrast, the block GMRES solvers used as a baseline, while effective on ideal Cartesian grids, show significant
sensitivity to geometric distortion, especially in three dimensions and for higher-order elements, where they
occasionally fail to converge for highly distorted meshes. While they can provide a simpler implementation, their
performance gap compared to the Braess--Sarazin approach widens as the difficulty of the problem increases. These
results confirm that the tight coupling of velocity and pressure within the Braess--Sarazin smoother is essential for
achieving robustness on the distorted or unstructured patches that comprise a general finite element mesh.

\subsection{Application in Geometric Multigrid}
Having established the performance of the local solvers on a single patch, we now return to their application as
smoothers within the global geometric multigrid algorithm. In a practical multigrid setting, the local problem on each
patch does not need to be solved to high accuracy. The goal of the smoother is to damp high-frequency error components,
which can typically be achieved with just a few iterations of an effective iterative method.

We denote by $N_{\text{MG}}$ the number of local iterations (e.g., p-multigrid V-cycles or Braess-Sarazin steps)
applied within each local patch solve. Based on the local solver analysis, we focus on the Braess-Sarazin smoother
configured with a single Richardson iteration for the inner Schur complement solve. This configuration offers the best
balance of cost and robustness. We compare it against a standard Block-Triangular (Block MG) preconditioner and a
direct (exact) solver on each patch.

We evaluate the performance on a 2D test problem. The implementation of the patch smoother in 2D is currently
matrix-based and unoptimized compared to the matrix-free 3D implementation, but it allows for rigorous testing of
numerical robustness. We start with a coarse grid of $3 \times 3$ cells. This specific discretization is chosen to
allow for a central cell that is isolated from the domain boundary, which we will use to test coefficient jumps. This
coarse grid is refined uniformly to generate a hierarchy of meshes (typically $L=4$ levels). To assess robustness
against grid quality, we introduce random geometric distortions. The interior vertices of the fine grid are perturbed
by a random vector of magnitude $\delta h$, where $h$ is the local cell width and $\delta \in \{0, 10\%, 25\%, 35\%\}$
is the distortion parameter.

\begin{figure}[htbp]
  \centering
  \begin{subfigure}[t]{0.32\textwidth}
    \centering
    \includegraphics[width=\textwidth]{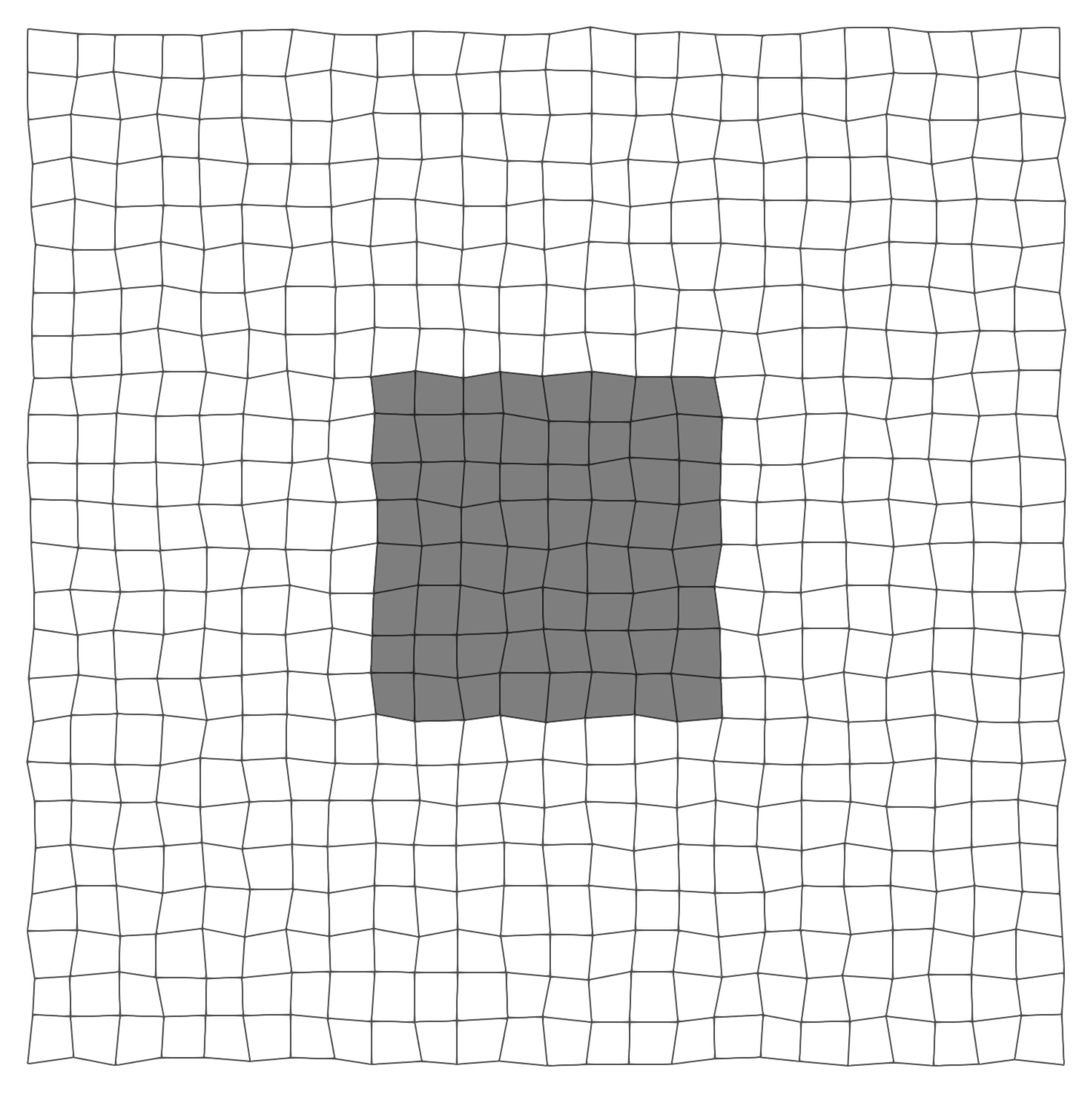}
    \caption{ $\delta= 10\%$.}
    \label{fig:distorted_10}
  \end{subfigure}
  \hfill
  \begin{subfigure}[t]{0.32\textwidth}
    \centering
    \includegraphics[width=\textwidth]{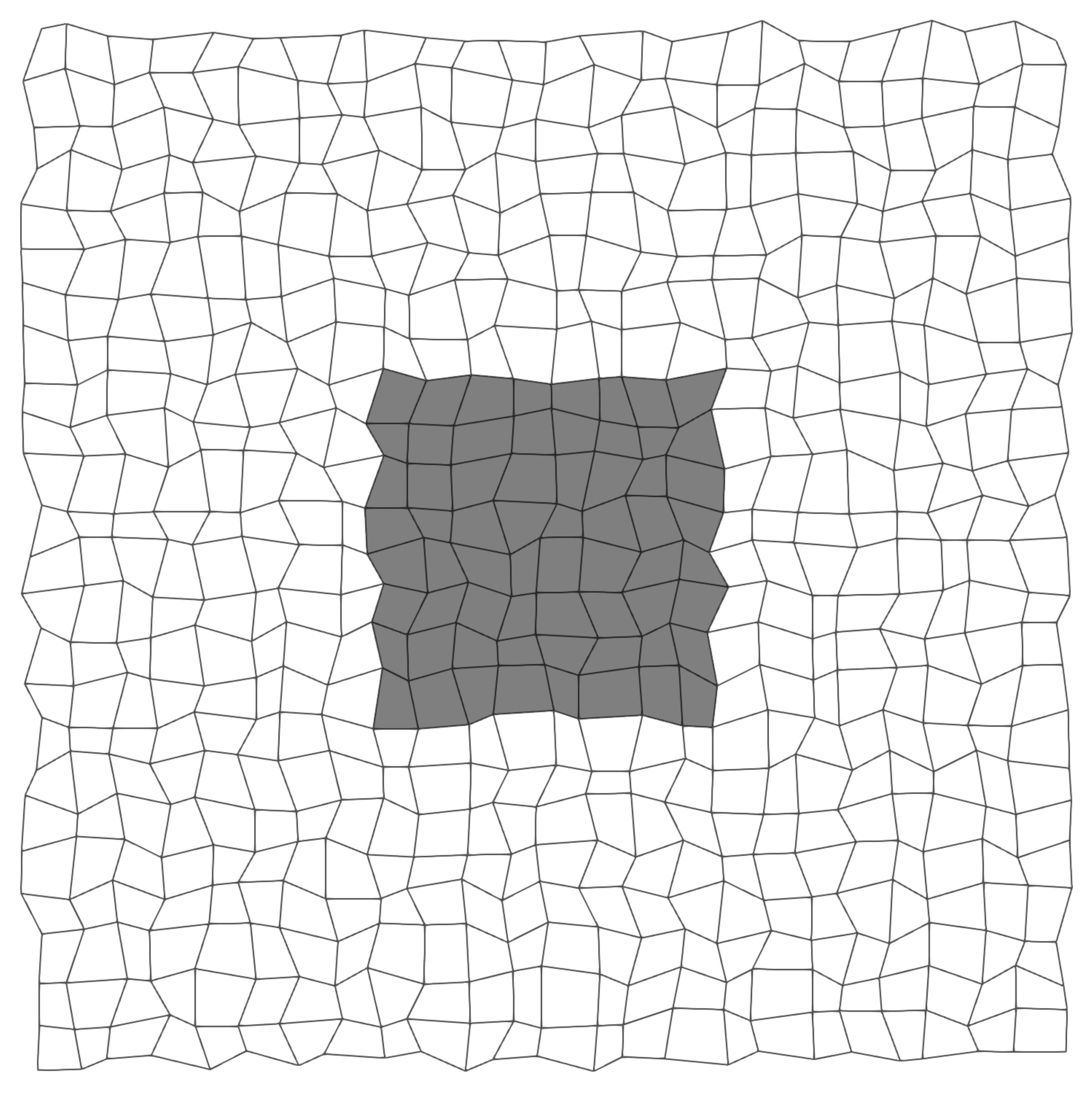}
    \caption{ $\delta= 25\%$.}
    \label{fig:distorted_25}
  \end{subfigure}
  \hfill
  \begin{subfigure}[t]{0.32\textwidth}
    \centering
    \includegraphics[width=\textwidth]{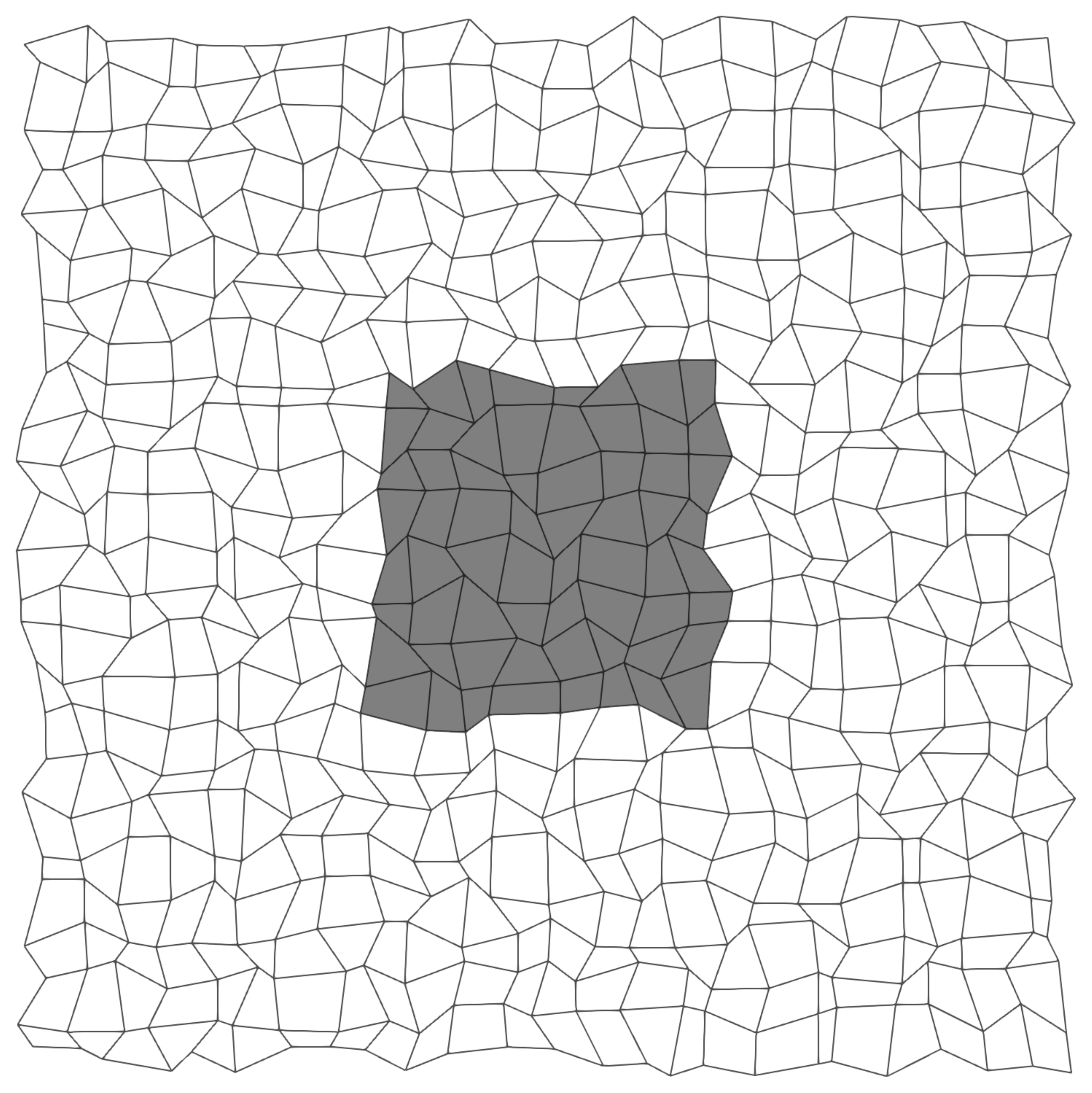}
    \caption{ $\delta= 35\%$.}
    \label{fig:distorted_35}
  \end{subfigure}
  \caption{Meshes with $L=4$ multigrid levels used in the numerical experiments to assess robustness to geometric distortion. 2D mesh distorted by various factors $\delta \in \{10\%, 25\%, 35\%\}$. The shaded region indicates the subdomain of high viscosity used in the coefficient-jump test.}
\end{figure}

The patch smoother is applied multiplicatively within a global geometric multigrid V-cycle with one pre- and one
post-smoothing step. We use FGMRES as the outer solver and iterate until the relative residual is reduced by $10^{-8}$.

Interestingly, the global GMRES iteration counts are notably low, often appearing even more favorable than the
convergence rates observed for isolated local problems. This effect is a direct consequence of the overlapping
structure of the vertex-patch smoother: in a multiplicative sweep over a structured $d$-dimensional mesh, each cell is
contained within $2^d$ different patches. As a result, the degrees of freedom on each cell are updated multiple times
during a single global smoothing step, significantly enhancing the error damping. A similar observation was made and
analyzed for the scalar Laplace problem in~\cite{wichrowski2025local}.

Table~\ref{tab:gmres_iterations_jacobi_2D} presents the results for the stationary Stokes problem with constant
coefficients. We observe that the simple Block MG approach is not robust. It fails to converge (NC) for
$N_{\text{MG}}=1$ and degrades significantly with mesh distortion and increasing polynomial degree $p$, even with more
local iterations. In sharp contrast, the Braess-Sarazin smoother with just a single step ($N_{\text{MG}}=1$) delivers
excellent performance. Its iteration counts are almost identical to those obtained using an exact local solver,
showcasing its ability to effectively handle the saddle-point nature of the equations locally. This robustness holds
true even for high polynomial degrees ($p=7$) and severe mesh distortion ($\delta=35\%$).

\begin{table}[htbp]
  \centering
  \caption{GMRES iteration counts (2D, problem with $L=4$ multigrid levels) for a geometric multigrid preconditioner
    using the patch smoother with a local p-multigrid solver. Columns show increasing mesh distortion $\delta$.
    The third column lists $N_{\text{MG}}$ (number of local p-multigrid cycles) for each row. Entries are FGMRES iteration counts to reduce the residual by $10^{-8}$. ``NC'' indicates no
    convergence within the iteration limit used in the experiments.}
  \begin{tabular}{l l c | c c c c}
    \toprule
    $p$ & Local solver                     & $N_{\text{MG}}$ & $\delta=0\%$ & $\delta=10\%$ & $\delta=25\%$ & $\delta=35\%$ \\
    \midrule
    \multirow{4}{*}{2}
        & \multirow{2}{*}{Block MG}        & 1               & NC           & NC            & NC            & NC            \\
        &                                  & 2               & 7            & 7             & 9             & 10            \\
    \arrayrulecolor{gray!80}\cline{2-7}\arrayrulecolor{black}
        & Braess--Sarazin                  & 1               & 7            & 7             & 8             & 9             \\
    \rowcolor{gray!15}
        & Exact                            & --              & 7            & 7             & 8             & 9             \\
    \midrule
    \multirow{5}{*}{3}
        & \multirow{3}{*}{Block MG}        & 1               & NC           & NC            & NC            & NC            \\
        &                                  & 2               & 7            & 7             & 9             & 12            \\
        &                                  & 3               & 6            & 7             & 8             & 17            \\
    \arrayrulecolor{gray!80}\cline{2-7}\arrayrulecolor{black}
        & Braess--Sarazin                  & 1               & 6            & 7             & 8             & 10            \\
    \rowcolor{gray!15}
        & Exact                            & --              & 6            & 7             & 8             & 9             \\
    \midrule
    \multirow{5}{*}{4}
        & \multirow{3}{*}{Block MG}        & 1               & NC           & NC            & NC            & NC            \\
        &                                  & 2               & 6            & 7             & 12            & 17            \\
        &                                  & 3               & 6            & 7             & 10            & NC            \\
    \arrayrulecolor{gray!80}\cline{2-7}\arrayrulecolor{black}
        & Braess--Sarazin                  & 1               & 6            & 7             & 8             & 9             \\
    \rowcolor{gray!15}
        & Exact                            & --              & 6            & 7             & 8             & 8             \\
    \midrule
    \multirow{6}{*}{7}
        & \multirow{3}{*}{Block MG}        & 1               & NC           & NC            & NC            & NC            \\
        &                                  & 2               & 12           & 22            & NC            & NC            \\
        &                                  & 3               & 11           & NC            & NC            & NC            \\
    \arrayrulecolor{gray!80}\cline{2-7}\arrayrulecolor{black}
        & \multirow{2}{*}{Braess--Sarazin} & 1               & 8            & 8             & 10            & 12            \\
        &                                  & 2               & 6            & 7             & 8             & 9             \\
    \rowcolor{gray!15}
        & Exact                            & --              & 5            & 6             & 7             & 8             \\
    \bottomrule
  \end{tabular}
  \label{tab:gmres_iterations_jacobi_2D}
\end{table}

A significant challenge for Stokes solvers is robustness in the presence of large jumps in material coefficients. To
test this, we utilize the central cell of our $3 \times 3$ coarse grid. We assign a viscosity of $\mu = 10^6$ to this
central cell (and all its fine-grid descendants), while the rest of the domain has $\mu = 1$.

Table~\ref{tab:gmres_iterations_mu1e6_2D} compares the performance of the Braess-Sarazin smoother against the exact
local solver in this setting. The results are compelling: the Braess-Sarazin smoother ($N_{\text{MG}}=1$) maintains
iteration counts that are remarkably close to the exact local solve, even with the $10^6$ jump in viscosity and large
mesh distortions. This confirms that the Braess-Sarazin approach correctly captures the local scalings induced by the
coefficient jump and separates them from the global smoothing process.

\begin{table}[htbp]
  \centering
  \caption{GMRES iteration counts (2D, $\mu=10^{6}$, with $L=4$ multigrid levels) comparing the Braess--Sarazin local
  smoother with one and two local p-multigrid cycles and the exact local solve. Columns show increasing
  mesh distortion $\delta$. Entries are FGMRES iteration counts to reduce the residual by $10^{-8}$. }
  \begin{tabular}{l l c | c c c c}
    \toprule
    $p$ & Local solver                     & $N_{\text{MG}}$ & $\delta=0\%$ & $\delta=10\%$ & $\delta=25\%$ & $\delta=35\%$ \\
    \midrule
    \multirow{3}{*}{2}
        & \multirow{2}{*}{Braess--Sarazin} & 1               & 6            & 7             & 9             & 9             \\
        &                                  & 2               & 6            & 7             & 9             & 10            \\
    \rowcolor{gray!15}
        & Exact                            & --              & 6            & 7             & 8             & 10            \\
    \midrule
    \multirow{3}{*}{3}
        & \multirow{2}{*}{Braess--Sarazin} & 1               & 6            & 7             & 8             & 10            \\
        &                                  & 2               & 6            & 7             & 8             & 9             \\
    \rowcolor{gray!15}
        & Exact                            & --              & 6            & 7             & 8             & 9             \\
    \midrule
    \multirow{3}{*}{4}
        & \multirow{2}{*}{Braess--Sarazin} & 1               & 6            & 6             & 8             & 9             \\
        &                                  & 2               & 6            & 6             & 7             & 8             \\
    \rowcolor{gray!15}
        & Exact                            & --              & 6            & 6             & 8             & 8             \\
    \midrule
    \multirow{3}{*}{7}
        & \multirow{2}{*}{Braess--Sarazin} & 1               & 8            & 8             & 10            & 12            \\
        &                                  & 2               & 6            & 6             & 8             & 9             \\
    \rowcolor{gray!15}
        & Exact                            & --              & 5            & 6             & 7             & 7             \\
    \bottomrule
  \end{tabular}
  \label{tab:gmres_iterations_mu1e6_2D}
\end{table}

\section{Conclusion}
\label{sec:conclusion}

In this work, we investigated the feasibility of using fully iterative, multigrid-based solvers for the local patch
problems in Stokes smoothers. Our goal was to determine if an \emph{inexact} solver could maintain global robustness
while avoiding the memory and computational costs of direct dense factorizations, which become prohibitive for
high-order discretizations.

Among the strategies tested, the {Braess-Sarazin smoother} emerged as a particularly strong candidate for the local
solve. Even when approximating the velocity inverse with a simple $p$-multigrid cycle, it demonstrated resilience to
mesh distortion and high-contrast viscosity coefficients. These results suggest that maintaining the saddle-point
structure within the patch is more beneficial than striving for high precision with decoupled block solvers. Moreover,
the efficiency of the Braess-Sarazin approach within the global multigrid hierarchy is striking: a single local
application ($N_{\text{MG}}=1$) yields global iteration counts that are virtually identical to those obtained using an
exact direct solver on the patches. In contrast, standard block solvers frequently failed to ensure global convergence
under similar conditions, even when allowed multiple local iterations.

It is important to note that while the numerical experiments in this paper were conducted using sparse matrix data
structures to facilitate rapid prototyping, the algorithms are designed specifically for a matrix-free context. The
iteration counts and robustness observed here serve as a validation for future high-performance implementations relying
on sum-factorization. Looking forward, this recursive smoothing strategy appears especially promising for
$H(\text{div})$-conforming discretizations, such as Raviart-Thomas elements. Since patch smoothers are often
theoretically required for these spaces, having a lean, matrix-free local solver could make high-order flux-continuous
schemes significantly more attractive for large-scale Stokes and Navier-Stokes simulations.

\section*{Acknowledgments}

The generative AI (Gemini, ChatGPT, Claude) was used to assist in drafting and proofreading parts of this manuscript.
Any remaining errors are the author's responsibility.

\noindent The author warmly thanks his mother, Grażyna, and his brother, Wiktor, for their support. Special thanks to Micro and Konda, the family cats, for their comforting companionship during the writing of this paper.

\bibliographystyle{siamplain}
\bibliography{literature}

@article{ArnoldFalkWinther00,
    AUTHOR = {Arnold, Douglas N. and Falk, Richard S. and Winther, Ragnar},
     TITLE = {Multigrid in {$H({\rm div})$} and {$H({\rm curl})$}},
   JOURNAL = {Numer. Math.},
  FJOURNAL = {Numer. Math.},
    VOLUME = 85,
      YEAR = 2000,
    NUMBER = 2,
     PAGES = {197--217},
       DOI = {10.1007/PL00005386}
}

@Book{Bramble93,
  title = {Multigrid Methods},
  author = {Bramble, J. H.},
  number = {294},
  publisher = {Longman Scientific},
  series = {Pitman research notes in mathematics series},
  year = {1993}
}

@article{burstedde2011p4est,
  title={p4est: Scalable algorithms for parallel adaptive mesh refinement on forests of octrees},
  author={Burstedde, Carsten and Wilcox, Lucas C and Ghattas, Omar},
  journal={SIAM J. Sci. Comput.},
  volume={33},
  number={3},
  pages={1103--1133},
  year={2011},
  publisher={SIAM}
}

@Article{dealii2019design,
        title   = {The {deal.II} finite element library: Design, features, and insights},
        author  = {Daniel Arndt and Wolfgang Bangerth and Denis Davydov and
                   Timo Heister and Luca Heltai and Martin Kronbichler and
                   Matthias Maier and Jean-Paul Pelteret and Bruno Turcksin and
                   David Wells},
        journal = {Comput. \& Math. Appl.},
        year    = {2021},
        DOI     = {10.1016/j.camwa.2020.02.022},
        pages   = {407-422},
        volume  = {81},
        url     = {https://arxiv.org/abs/1910.13247}
      }

@article{dealII97,
  title   = {The deal.{II} Library, Version 9.7 },
  author  = {Daniel Arndt and  Wolfgang Bangerth and Maximilian Bergbauer and 
   Bruno Blais and  Marc Fehling and  Rene Gassmöller 
   and  Timo Heister and  Luca Heltai and  Martin Kronbichler
    and  Matthias Maier and  Peter Munch and  Sam Scheuerman
     and  Bruno Turcksin and  Siarhei Uzunbajakau and  David Wells and  Michał Wichrowski},
  journal = {preprint},
  year    = {2025},
  url     = {https://dealii.org/deal97-preprint.pdf}
}

@article{brubeck2021scalable,
author = {Brubeck, Pablo D. and Farrell, Patrick E.},
title = {A Scalable and Robust Vertex-Star Relaxation for High-Order {FEM}},
journal = {SIAM J. Sci. Comput.},
volume = {44},
number = {5},
pages = {A2991-A3017},
year = {2022},
doi = {10.1137/21M1444187}
}

@Book{Hackbusch85,
  title = {Multi-grid Methods and Applications},
  author = {Hackbusch, W.},
  publisher = {Springer},
  address = {Heidelberg},
  year = 1985
}

@Article{JanssenKanschat11,
  author =       {Janssen, B. and Kanschat, G.},
  title =        {Adaptive multilevel methods with local smoothing for
                  {${H}^1$}- and {$H^{\text{curl}}$}-conforming high order
                  finite element methods},
  journal =      {SIAM J. Sci. Comput.},
  year =         2011,
  volume =       33,
  number =       4,
  pages =        {2095--2114},
  doi = {10.1137/090778523}
}

@Article{KanschatMao15,
  author = 	 {Kanschat, G. and Mao, Y.},
  title = 	 {Multigrid methods for {$\mathbf H^{\text{div}}$}-conforming discontinuous
{G}alerkin methods for the {S}tokes equations},
  journal = 	 {J. Numer. Math.},
  year = 	 2015,
  volume = 	 23,
  number = 	 1,
  pages = 	 {51--66},
  doi =	 {10.1515/jnma-2015-0005}
}

@article{Kronbichler2012,
author = {Kronbichler, Martin and Kormann, Katharina},
journal = {Computers \& Fluids},
pages = {135--147},
title = {{A generic interface for parallel cell-based finite element operator application}},
volume = {63},
year = {2012},
doi = {10.1016/j.compﬂuid.2012.04.012}
}

@article{Kronbichler2017a,
 author = {Kronbichler, Martin and Kormann, Katharina},
 title = {Fast Matrix-Free Evaluation of Discontinuous {G}alerkin Finite Element Operators},
 year = {2019},
 volume = {45},
 number = {3},
 journal = {ACM Trans. Math. Softw.},
 pages = {29/1--40},
 doi = {10.1145/3325864}
}

@article{Lynch1964,
author = {Lynch, R. E. and Rice, J. R. and Thomas, D. H.},
title = {Direct solution of partial difference equations by tensor product methods},
journal = {Numer. Math.},
volume = 6,
pages = {185--199},
year = 1964,
doi = {10.1007/BF01386067}
}

@Article{WitteArndtKanschat21,
  author =       {Witte, J. and Arndt, D. and Kanschat, G.},
  title =        {Fast Tensor Product {S}chwarz Smoothers for High-Order Discontinuous {G}alerkin Methods},
  journal =      {Comput. Meth. Appl. Math.},
  number = {3},
  volume = {21},
  year = {2021},
  pages = {709--728},
  doi = {10.1515/cmam-2020-0078}
}

@article{kronbichler2019multigrid,
  title={Multigrid for matrix-free high-order finite element computations on graphics processors},
  author={Kronbichler, Martin and Ljungkvist, Karl},
  journal={ACM Trans. Parallel Comput.},
  volume={6},
  number={1},
  pages={2/1--32},
  year={2019},
  doi={10.1145/3322813}
}

@article{thomas2003textbook,
  title={Textbook multigrid efficiency for fluid simulations},
  author={Thomas, James L and Diskin, Boris and Brandt, Achi},
  journal={Annu. Rev. Fluid Mech.},
  volume={35},
  number={1},
  pages={317--340},
  year={2003},
  doi={10.1146/annurev.fluid.35.101101.161209}
}

@article{hong2016robust,
  title={A robust multigrid method for discontinuous {G}alerkin discretizations of {S}tokes and linear elasticity equations},
  author={Hong, Qingguo and Kraus, Johannes and Xu, Jinchao and Zikatanov, Ludmil},
  journal={Numer. Math.},
  volume={132},
  number={1},
  pages={23--49},
  year={2016},
  publisher={Springer}
}

@article{munch2023cache,
 author = {Peter Munch and Martin Kronbichler},
 title = {Cache-optimized and low-overhead implementations of additive {S}chwarz methods for high-order {FEM} multigrid computations},
 journal = {Int. J. High Perf. Comput. Appl.},
 year = 2023,
 note = {In press},
 doi = {10.1177/10943420231217221}
}

@article{wichrowski2025smoothers,
  title={Smoothers with Localized Residual Computations for Geometric Multigrid Methods for Higher-Order Finite Elements},
  author={Wichrowski, Micha{\l} and Munch, Peter and Kronbichler, Martin and Kanschat, Guido},
  journal={SIAM Journal on Scientific Computing},
  volume={47},
  number={3},
  pages={B645--B664},
  year={2025},
  publisher={SIAM}
}

@article{wichrowski2022matrix,
  title={A matrix-free multilevel preconditioner for the generalized {S}tokes problem with discontinuous viscosity},
  author={Wichrowski, Micha{\l} and Krzy{\.z}anowski, Piotr},
  journal={Journal of Computational Science},
  volume={63},
  pages={101804},
  year={2022},
  publisher={Elsevier}
}

@article{cui2025implementation,
  title={An implementation of tensor product patch smoothers on {GPU}s},
  author={Cui, Cu and Grosse-Bley, Paul and Kanschat, Guido and Strzodka, Robert},
  journal={SIAM Journal on Scientific Computing},
  volume={47},
  number={2},
  pages={B280--B307},
  year={2025},
  publisher={SIAM}
}

@article{witte2025tensor,
  title={Tensor-product vertex patch smoothers for biharmonic problems},
  author={Witte, Julius and Cui, Cu and Bonizzoni, Francesca and Kanschat, Guido},
  journal={Computational Methods in Applied Mathematics},
  number={0},
  year={2025},
  publisher={De Gruyter}
}

@article{pavarino1993additive,
  title={Additive {S}chwarz methods for the p-version finite element method},
  author={Pavarino, Luca F},
  journal={Numerische Mathematik},
  volume={66},
  number={1},
  pages={493--515},
  year={1993},
  publisher={Springer}
}

@article{pazner2020efficient,
  title={Efficient low-order refined preconditioners for high-order matrix-free continuous and discontinuous {G}alerkin methods},
  author={Pazner, Will},
  journal={SIAM Journal on Scientific Computing},
  volume={42},
  number={5},
  pages={A3055--A3083},
  year={2020},
  publisher={SIAM}
}

@article{margenberg2025hp,
  title={An hp multigrid approach for tensor-product space-time finite element discretizations of the {S}tokes equations},
  author={Margenberg, Nils and Bause, Markus and Munch, Peter},
  journal={arXiv preprint arXiv:2502.09159},
  year={2025}
}

@article{braess1997efficient,
  title={An efficient smoother for the {S}tokes problem},
  author={Braess, Dietrich and Sarazin, Regina},
  journal={Applied Numerical Mathematics},
  volume={23},
  number={1},
  pages={3--19},
  year={1997},
  publisher={Elsevier}
}

@article{jodlbauer2024matrix,
  title={Matrix-free monolithic multigrid methods for {S}tokes and generalized {S}tokes problems},
  author={Jodlbauer, Daniel and Langer, Ulrich and Wick, Thomas and Zulehner, Walter},
  journal={SIAM Journal on Scientific Computing},
  volume={46},
  number={3},
  pages={A1599--A1627},
  year={2024},
  publisher={SIAM}
}

@article{zulehner2000class,
  title={A class of smoothers for saddle point problems},
  author={Zulehner, Walter},
  journal={Computing},
  volume={65},
  number={3},
  pages={227--246},
  year={2000},
  publisher={Springer-Verlag Wien}
}

@article{wichrowski2025geometric,
  title={A Geometric Multigrid Preconditioner for Discontinuous {G}alerkin Shifted Boundary Method},
  author={Wichrowski, Michal},
  journal={arXiv preprint arXiv:2506.12899},
  year={2025}
}

@article{cui2025multigrid,
  title={A multigrid method for CutFEM and its implementation on {GPU}},
  author={Cui, Cu and Kanschat, Guido},
  journal={arXiv preprint arXiv:2508.11608},
  year={2025}
}

@article{bergbauer2025high,
  title={High-performance matrix-free unfitted finite element operator evaluation},
  author={Bergbauer, Maximilian and Munch, Peter and Wall, Wolfgang A and Kronbichler, Martin},
  journal={SIAM Journal on Scientific Computing},
  volume={47},
  number={3},
  pages={B665--B689},
  year={2025},
  publisher={SIAM}
}

@article{wichrowski2025local,
  title={Local Solvers for High-Order Patch Smoothers via p-Multigrid},
  author={Wichrowski, Micha{\l}},
  journal={arXiv preprint arXiv:2510.17785},
  year={2025}
}

@book{deville2002high,
  title={High-order methods for incompressible fluid flow},
  author={Deville, Michel O and Fischer, Paul F and Mund, Ernest H},
  volume={9},
  year={2002},
  publisher={Cambridge university press}
}

@article{voronin2025monolithic,
  title={Monolithic algebraic multigrid preconditioners for the {S}tokes equations},
  author={Voronin, Alexey and MacLachlan, Scott and Olson, Luke N and Tuminaro, Raymond S},
  journal={SIAM Journal on Scientific Computing},
  volume={47},
  number={1},
  pages={A343--A373},
  year={2025},
  publisher={SIAM}
}

@article{olshanskii2006uniform,
  title={Uniform preconditioners for a parameter dependent saddle point problem with application to generalized {S}tokes interface equations},
  author={Olshanskii, Maxim A and Peters, J{\"o}rg and Reusken, Arnold},
  journal={Numerische Mathematik},
  volume={105},
  number={1},
  pages={159--191},
  year={2006},
  publisher={Springer}
}

@article{wichrowski2025pMGimplementaion,
  title={Multigrid p-Robustness at Jacobi Speeds: Efficient Matrix-Free Implementation of Local p-Multigrid Solvers},
  author={Wichrowski, Micha{\l}},
  journal={arXiv preprint arXiv:2512.02577},
  year={2025}
}

@article{cui2025Stokes,
  title={Multigrid methods for the {S}tokes problem on {GPU} systems},
  author={Cui, Cu and Kanschat, Guido},
  journal={Computers \& Fluids},
  pages={106703},
  year={2025},
  publisher={Elsevier}
}

@article{wichrowski2023exploiting,
  title={Exploiting high-contrast {S}tokes preconditioners to efficiently solve incompressible fluid--structure interaction problems},
  author={Wichrowski, Micha{\l} and Krzy{\.z}anowski, Piotr and Heltai, Luca and Stupkiewicz, Stanis{\l}aw},
  journal={International Journal for Numerical Methods in Engineering},
  volume={124},
  number={24},
  pages={5446--5470},
  year={2023},
  publisher={Wiley Online Library}
}

@article{chen2015multigrid,
  title={Multigrid methods for saddle point systems using constrained smoothers},
  author={Chen, Long},
  journal={Computers \& Mathematics with Applications},
  volume={70},
  number={12},
  pages={2854--2866},
  year={2015},
  publisher={Pergamon}
}

@article{krzyzanowski2011block,
  title={On block preconditioners for saddle point problems with singular or indefinite (1, 1) block},
  author={Krzy{\.z}anowski, Piotr},
  journal={Numerical Linear Algebra with Applications},
  volume={18},
  number={1},
  pages={123--140},
  year={2011},
  publisher={Wiley Online Library}
}

@article{krzyzanowski2001block,
  title={On block preconditioners for nonsymmetric saddle point problems},
  author={Krzyzanowski, Piotr},
  journal={SIAM Journal on Scientific Computing},
  volume={23},
  number={1},
  pages={157--169},
  year={2001},
  publisher={SIAM}
}

@article{rudi2017weighted,
  title={Weighted {BFBT} preconditioner for {S}tokes flow problems with highly heterogeneous viscosity},
  author={Rudi, Johann and Stadler, Georg and Ghattas, Omar},
  journal={SIAM Journal on Scientific Computing},
  volume={39},
  number={5},
  pages={S272--S297},
  year={2017},
  publisher={Society for Industrial and Applied Mathematics}
}

@article{vanka1986block,
  title={Block-implicit multigrid solution of {N}avier-{S}tokes equations in primitive variables},
  author={Vanka, S Pratap},
  journal={Journal of Computational Physics},
  volume={65},
  number={1},
  pages={138--158},
  year={1986},
  publisher={Academic Press}
}

@article{borzacchiello2017box,
  title={Box-relaxation based multigrid solvers for the variable viscosity {S}tokes problem},
  author={Borzacchiello, Domenico and Leriche, Emmanuel and Blotti{\`e}re, Beno{\^\i}t and Guillet, Jacques},
  journal={Computers \& Fluids},
  volume={156},
  pages={515--525},
  year={2017},
  publisher={Pergamon}
}

@article{may2015scalable,
  title={A scalable, matrix-free multigrid preconditioner for finite element discretizations of heterogeneous {S}tokes flow},
  author={May, Dave A and Brown, Jed and Le Pourhiet, Laetitia},
  journal={Computer methods in applied mechanics and engineering},
  volume={290},
  pages={496--523},
  year={2015},
  publisher={North-Holland}
}

@article{wobker2009numerical,
  title={Numerical studies of {V}anka-type smoothers in computational solid mechanics},
  author={Wobker, Hilmar and Turek, Stefan},
  journal={Advances in Applied Mathematics and Mechanics},
  volume={1},
  number={1},
  pages={29--55},
  year={2009}
}

@article{farrell2021local,
  title={A local {F}ourier analysis of additive {V}anka relaxation for the {S}tokes equations},
  author={Farrell, Patrick E and He, Yunhui and MacLachlan, Scott P},
  journal={Numerical Linear Algebra with Applications},
  volume={28},
  number={3},
  pages={e2306},
  year={2021},
  publisher={Wiley Online Library}
}

@article{harper2023compression,
  title={Compression and Reduced Representation Techniques for Patch-Based Relaxation},
  author={Harper, Graham and Tuminaro, Ray},
  journal={arXiv preprint arXiv:2306.10025},
  year={2023}
}

@article{margenberg2025multigrid,
  title={An Multigrid Approach for Tensor-Product Space-Time Finite Element Discretizations of the {S}tokes Equations},
  author={Margenberg, Nils and Bause, Markus and Munch, Peter},
  journal={SIAM Journal on Scientific Computing},
  volume={47},
  number={6},
  pages={B1503--B1529},
  year={2025},
  publisher={SIAM}
}

@incollection{RaviartThomas,
    AUTHOR = {Raviart, Pierre-Arnaud and Thomas, Jean-Marie},
    EDITOR = {Galligani, Ilio and Magenes, Enrico},
     TITLE = {A mixed finite element method for 2nd order elliptic problems},
 BOOKTITLE = {Mathematical aspects of finite element methods},
    VOLUME = {606},
      YEAR = {1977},
       DOI = {10.1007/BFb0064470},
     PAGES = {{292--315}},
}

@article{RaviartThomas2,
    AUTHOR = {N\'ed\'elec, Jean-Claude},
     TITLE = {Mixed finite elements in \(\mathbb{R}^3\)},
   JOURNAL = {Numerische Mathematik},
    VOLUME = {35},
    NUMBER = {3},
      YEAR = {1980},
       DOI = {10.1007/BF01396415},
     PAGES = {{315--341}},
}

\newpage
\appendix
\section{Closed form of the local inverse approximate and computational complexity estimate}
\label{sec:appendix-local-inverse}

The local solver via the $p$-multigrid method employing a Braess--Sarazin smoother on each level yields an approximate
inverse of the Stokes operator on a single patch that can be written down explicitly. In this appendix, we derive the
closed-form expression for this operator to illustrate that every component can be evaluated solely using Kronecker
products of one-dimensional matrices, and also for a bit of mathematical fun.

\subsection{Level and Transfer Operators}

Let $k = 1, \dots, L$ denote the local p-multigrid levels, corresponding to polynomial degrees $p_1 < p_2 < \dots <
    p_L$. We denote the finite element spaces on level $k$ as $\mathbb{V}_k \times \mathbb{Q}_k$. Note that for the
coarsest level $k=1$ (corresponding to $p=1$), the pressure space is effectively empty due to the mean value constraint
on a single macro-element, so we perform calculations only on the velocity space $\mathbb{V}_1$.

The discrete Stokes operator on level $k$ is given by
\begin{equation}
    \mathcal{A}_k = \begin{pmatrix} A_k & B_k^T \\ B_k & 0 \end{pmatrix}.
\end{equation}
For $k=1$, this reduces to $\mathcal{A}_1 = A_1$.

We define the prolongation operators $\mathcal{P}_{k-1}^k: (\mathbb{V}_{k-1} \times \mathbb{Q}_{k-1}) \to (\mathbb{V}_k
    \times \mathbb{Q}_k)$ as the natural finite element embeddings. We assume that the finite element spaces are
constructed as tensor products of one-dimensional spaces. This property is satisfied by common choices such as the
$\mathbb{Q}_p - \mathbb{Q}^{DG}_{p-2}$ pair as well as Raviart--Thomas of degree $p$ elements and their corresponding
pressure spaces $\mathbb{Q}^{DG}_{p}$.

The mentioned operators exhibit a tensor-product structure. Specifically, for a $d$-dimensional patch, the transfer
operator can be written as a Kronecker product of 1D transfer matrices:
\begin{equation}
    P_{k-1}^k = P_{1D} \otimes \dots \otimes P_{1D},
\end{equation}
where $P_{1D} \in \mathbb{R}^{(p_{k}+1) \times (p_{k-1}+1)}$ is the 1D embedding matrix. The restriction operators $\mathcal{R}_k^{k-1}$ are the adjoints of the prolongations.

The efficient implementation of the local solver should rely on the fact that both transfer operators and level
operators can be evaluated using sum factorization. To define the level operators $A_k$ and $B_k$, let $V_k, D_k \in
    \mathbb{R}^{n_q \times (p_k+1)}$ be the 1D Vandermonde (evaluation) and derivative matrices at the quadrature points,
and let $V_{p,k} \in \mathbb{R}^{n_q \times p_k}$ be the corresponding evaluation matrix for the pressure space
($\mathbb{Q}_{p-1}$). In 3D, we define the gradient operators on the reference cell as:
\begin{equation}
    \hat{\nabla}_1 = D_k \otimes V_k \otimes V_k, \quad \hat{\nabla}_2 = V_k \otimes D_k \otimes V_k, \quad \hat{\nabla}_3 = V_k \otimes V_k \otimes D_k.
\end{equation}
Let $J = \frac{\partial x}{\partial \xi}$ be the mapping Jacobian, $G = J^{-T}$ be its inverse transpose, and $W =
    \text{diag}(w_q) \det(J)$ be the diagonal matrix of quadrature weights scaled by the determinant of the Jacobian.

The physical gradients are then given by $\nabla_j = \sum_{i=1}^3 G_{ij} \hat{\nabla}_i$. The velocity block $A_{k,1}$
(for a single component) and the divergence blocks $B_{k,j}$ can be expressed as:
\begin{equation}
    A_{k,1} = \sum_{j=1}^3 \nabla_j^T W \nabla_j, \quad B_{k,j} = \mathbf{V}_{p,k}^T W \nabla_j,
\end{equation}
where $\mathbf{V}_{p,k} = V_{p,k} \otimes V_{p,k} \otimes V_{p,k}$ is the 3D pressure evaluation operator.

Since $W$ and $G$ are constant, the evaluation remains a structured sequence of 1D operations. For instance, the action
of one reference gradient component on a vector $u$ is evaluated as:
\begin{equation}
    (V_k \otimes V_k \otimes D_k) u = (V_k \otimes I \otimes I) (I \otimes V_k \otimes I) (I \otimes I \otimes D_k) u.
\end{equation}
This decomposition reduces the computational complexity from $O((p+1)^{6})$ to $O((p+1)^{4})$, which is crucial for
high polynomial degrees. Note that these estimates require the
underlying elements to have a tensor-product structure. This includes Raviart--Thomas elements, which, despite their more
complex definition, maintain this structure and thus benefit from the same gains in efficiency\cite{cui2025Stokes}.

\subsection{The Braess--Sarazin Smoother}

On each level $k > 1$, we apply a smoother $\mathcal{S}_k$. The action of the inverse $\mathcal{S}_k^{-1}$ corresponds
to one iteration of the Braess--Sarazin method with approximated Schur complement. Let $r = (r_u, r_p)^T$ be the
residual. The smoothing step $x_{new} = x_{old} + \mathcal{S}_k^{-1}(r - \mathcal{A}_k x_{old})$ (with initial guess
zero for the correction) produces a correction $(\delta u, \delta p)^T$ computed as follows:
\begin{enumerate}
    \item $\delta u^{(1)} = \tilde{A}_k^{-1} r_u$,
    \item $\delta p = \tilde{S}_k^{-1} (B_k \delta u^{(1)} - r_p)$,
    \item $\delta u = \delta u^{(1)} - \tilde{A}_k^{-1} B_k^T \delta p$.
\end{enumerate}
Here, $\tilde{A}_k$ and $\tilde{S}_k$ are diagonal scaling matrices incorporating relaxation parameters $\omega_a$ and $\omega_s$:
\begin{equation}
    \tilde{A}_k^{-1} = \omega_a (\text{diag}(A_k))^{-1}, \quad \tilde{S}_k^{-1} = \omega_s (\text{diag}(S_k))^{-1},
\end{equation}
where $S_k = B_k \tilde{A}_k^{-1} B_k^T$ is the Schur complement associated with the preconditioner $\tilde{A}_k$.

We can express the action of $\mathcal{S}_k^{-1}$ in block matrix form. Substituting step 1 into 2, we get
\begin{equation}
    \delta p = \tilde{S}_k^{-1} B_k \tilde{A}_k^{-1} r_u - \tilde{S}_k^{-1} r_p.
\end{equation}
Then, substituting $\delta p$ into step 3:
\begin{align}
    \delta u & = \tilde{A}_k^{-1} r_u - \tilde{A}_k^{-1} B_k^T (\tilde{S}_k^{-1} B_k \tilde{A}_k^{-1} r_u - \tilde{S}_k^{-1} r_p)                     \\
             & = (\tilde{A}_k^{-1} - \tilde{A}_k^{-1} B_k^T \tilde{S}_k^{-1} B_k \tilde{A}_k^{-1}) r_u + \tilde{A}_k^{-1} B_k^T \tilde{S}_k^{-1} r_p.
\end{align}
Thus, the operator $\mathcal{S}_k^{-1}$ is explicitly given by
\begin{equation}
    \mathcal{S}_k^{-1} = \begin{pmatrix}
        \tilde{A}_k^{-1} - \tilde{A}_k^{-1} B_k^T \tilde{S}_k^{-1} B_k \tilde{A}_k^{-1} & \tilde{A}_k^{-1} B_k^T \tilde{S}_k^{-1} \\
        \tilde{S}_k^{-1} B_k \tilde{A}_k^{-1}                                           & - \tilde{S}_k^{-1}
    \end{pmatrix}.
\end{equation}

\subsection{The p-Multigrid Operator}

Let $\mathcal{M}_k^{-1}$ denote the approximate inverse of $\mathcal{A}_k$ provided by the $p$-multigrid V-cycle with
one pre-smoothing and one post-smoothing step per level.

On the coarsest level $k=1$, the problem is solved exactly for the velocity:
\begin{equation}
    \mathcal{M}_1^{-1} = A_1^{-1}.
\end{equation}
Note that $A_1$ is a $3 \times 3$ matrix in 3D and $2 \times 2$ matrix in 2D, corresponding to the velocity space spanned on a single node.

For $k > 1$, the V-cycle operator is defined recursively. The inverse approximate consists of pre-smoothing, followed
by the coarse grid correction, and finally post-smoothing:
\begin{equation}
    \mathcal{M}_k^{-1} = (2 I - \mathcal{S}_k^{-1} \mathcal{A}_k) \mathcal{S}_k^{-1} + (I - \mathcal{S}_k^{-1} \mathcal{A}_k) \mathcal{P}_{k-1}^k \mathcal{M}_{k-1}^{-1} \mathcal{R}_k^{k-1} (I - \mathcal{A}_k \mathcal{S}_k^{-1}).
\end{equation}

For polynomial degrees $k=2,3$ the above formula does not require any further recursion since the coarsest level $k=1$
is reached, where the exact inverse is applied. Hence, we have obtained a closed-form expression for
$\mathcal{M}_k^{-1}$ for those cases, while for higher polynomial degrees the recursion should be expanded further.

\subsection{Complexity Analysis}

The computational complexity of the $p$-multigrid cycle is characterized by the number of applications of the discrete
velocity Laplacian $A_k$ and the divergence operator $B_k$ (and its adjoint $B_k^T$) on the finest level. We assume
that applying the diagonal inverse of $A_k$ (and similarly for the Schur complement) is negligible compared to the cost
of the matrix-vector multiplications. Consequently, we quantify the complexity by the total number $N_A$ and $N_B$ of
applications of $A_k$ and $B_k$, respectively. We further assume that the cost of applying $B_k$ and $B_k^T$ is the
same.

When a single Richardson iteration ($n_S=1$) is used for the inner Schur complement solve, the application of
$\mathcal{S}_k^{-1}$ involves exactly one multiplication by $B_k$ to form the pressure right-hand side and one by
$B_k^T$ to update the velocity correction. Thus, one full smoothing step consisting of an operator application
$\mathcal{A}_k$ followed by the application of $\mathcal{S}_k^{-1}$ requires $1 A_k$, $2 B_k$, and $2 B_k^T$.

Consider a V-cycle with $m$ pre-smoothing and $m$ post-smoothing steps. In addition to the $2m$ smoothing steps, the
multigrid cycle requires one additional application of $\mathcal{A}_k$ to compute the defect for the restriction to the
coarser levels. The total number of multiplications on the finest level is:
\begin{equation}
    N_A = 2m + 1, \quad N_B = 2m + (2m + 1) = 4m + 1.
\end{equation}
For the standard configuration used throughout this work ($m=1$, $n_S=1$), the cost on the finest level is exactly
$N_A=3$ and $N_B=5$.

We also express the cost in terms of sum-factorization kernels, where each kernel represents the application of a
$d$-dimensional tensor-product operator (e.g., $X \otimes Y \otimes Z$ in 3D). For a stationary Stokes problem,
evaluation of $ A_k$ requires 18 sum-factorization kernels (9 gradient components times 2 for the application and its
transpose), while evaluation of $ B_k$ requires 10 sum-factorization kernels (9 gradient components to compute
divergence and one to multiply by pressure and integrate). This results in the cost of operator $\mathcal{A}_k$ being
38 sum-factorization kernels while the cost of all smoothing steps (one pre- and one post-smoothing step) is 116
sum-factorization kernels. The cost of transfer operators requires an additional 6 sum-factorization kernels (3 for
prolongation and 3 for restriction) while the cost of evaluation on the coarser levels operators is less then
\nicefrac{1}{4} of the cost of the finest level. Thus, the total cost of the $p$-multigrid V-cycle with one pre- and
one post-smoothing step cost approximately the same as 152 sum-factorization kernels on the finest level, that is
around 4.0 times the cost of a single application of the discrete Stokes operator $\mathcal{A}_k$.
In~\cite{wichrowski2025pMGimplementaion}, empirical measurements for the Poisson equation showed a cost of
approximately 2.27 times the operator application.

Note that mileage may vary: on modern architectures some operations are cheaper than others so the cost will be
hardware dependent. Additionally, in the application of $\mathcal{A}_k$ some operations may be fused, or lower
precision computations may be used in p-multigrid cycle to reduce the effective cost, so this is a rough estimate.
Finally, on Cartesian grids some optimizations are possible.

\end{document}